\documentclass{elsarticle}
\usepackage{amsmath,amssymb,amsfonts,amsthm,mathrsfs}

\usepackage{graphicx,color}

\theoremstyle{plain}
\numberwithin{equation}{section}

\newtheorem{theorem}{Theorem}[section]
\newtheorem{proposition}[theorem]{Proposition}
\newtheorem{lemma}[theorem]{Lemma}
\newtheorem{corollary}[theorem]{Corollary}
\newtheorem{remark}[theorem]{Remark}

\def\im{{\rm Im\, }}
\def\codim{{\rm codim\,}}
\def\kr{{\rm Ker\, }}
\def\diag{{\rm diag\, }}

\def\rank{{\rm rank\, }}

\def\ind{{\rm ind\,}}

\def\BC{{\mathbb C}}
\def\BD{{\mathbb D}}

\definecolor{green}{rgb}{.0,.6,.0}

\newcommand{\bpr}{{\noindent\textbf{Proof.}\ \ }}
\newcommand{\epr}{{\hfill $\Box$}}

\begin{document}

\begin{frontmatter}
\title {Wiener-Hopf indices of unimodular functions on the unit circle, revisited.}

\author[AF]{A.E. Frazho}\ead{frazho@purdue.edu}
\author[MAK]{M.A. Kaashoek}
\author[AR]{A.C.M. Ran}\ead{a.c.m.ran@vu.nl}
\author[FvS]{F. van Schagen}\ead{f.van.schagen@vu.nl}

\address[AF]{Department of Aeronautics and Astronautics,
Purdue University,
West Lafayette, IN 47907, USA.}

\address[MAK]{Department of Mathematics, Faculty of Science, VU Amsterdam,
De Boelelaan 1111, 1081 HV Amsterdam}

\address[AR]{Department of Mathematics, Faculty of Science, VU Amsterdam, De Boelelaan 1111, 1081 HV Amsterdam, The Netherlands and Research Focus: Pure and Applied Analytics, North-West~University, Potchefstroom, South Africa}

\address[FvS]{Department of Mathematics, Faculty of Science, VU Amsterdam,
De Boelelaan 1111, 1081 HV Amsterdam}

\begin{abstract}
Inspired by the  paper of Groenewald, Kaashoek and Ran
(Wiener-Hopf indices of unitary functions on the unit circle in terms of realizations and
related results on Toeplitz operators. \emph{Indag. Math.} 28, (2017),  649-710),
we present an operator-theoretic approach to provide further insight
and simpler  computational formulas  for the Wiener-Hopf indices of a rational matrix valued function
 taking unimodular values on the unit circle.
\end{abstract}


\end{frontmatter}

{\date{}}

\begin{center}
\emph{Dedicated to Albrecht B\"{o}ttcher on the occasion of his seventieth  birthday.}
\end{center}


\setcounter{equation}{0}
\section{Introduction}

Wiener-Hopf factorization of matrix valued functions plays an important role in determining the Fredholm properties of several classes of operators, such as singular integral operators and Toeplitz operators, see e.g., \cite{BS, CGOT3,GF,GGKOT49,GGKOT63} and \cite{GKS}. To make this more explicit, let $R(z)$ be an $m\times m$ matrix valued function on the unit circle, which is
continuous and takes invertible values for $z$ in the unit circle $\mathbb{T}$. A factorization
\[
R(z)=W_-(z)\diag (z^{\kappa_j})_{j=1}^m W_+(z),
\]
where $W_-$ and its inverse are analytic outside the open unit disc, including infinity, and $W_+$ and its inverse are analytic inside the closed unit disc, and $\kappa_j\in\mathbb{Z}$ for $j=1, \ldots, m$,
is called a (\emph{right}) \emph{Wiener-Hopf factorization} with respect to the unit circle.
The integers $\kappa_j$ are uniquely determined by $R$ and they are called the \emph{Wiener-Hopf indices} of $R$.
The Toeplitz operator $T_R$ with symbol $R$ is Fredholm, and it factorizes accordingly as $T_R=T_{W_-}T_DT_{W_+}$, where $D(z)=\diag(z^{\kappa_j})_{j=1}^m$. Moreover, $T_{W_-}$ and $T_{W_+}$ are invertible operators. In that case, the dimension of the kernel of $T_R$ is equal to the dimension of the kernel of $T_D$, and the codimension of the range of $T_R$ is equal to the codimension of the range of $T_D$. In fact, these dimensions are completely determined by the Wiener-Hopf indices $\kappa_j$ as follows:
\[
\mathfrak{n}(T_R):=\dim(\kr(T_R))=\sum_{\kappa_j\leq -1} -\kappa_j, \quad d(T_R):=\codim(\im (T_R))=\sum_{\kappa_j\geq 1} \kappa_j.
\]
Conversely, the Wiener-Hopf indices $\kappa_j$ can be computed from dimensions of kernels and dimensions of cokernels as follows.

Let $ -\kappa_1 , -\kappa_2  , \ldots, -\kappa_p  $ with
$ \kappa_1 \geq \kappa_2 \geq \cdots \geq \kappa_p  $ be the negative
Wiener-Hopf indices of the function $ R $.
Then (see \cite{GGKOT63}  Theorem XXIV.4.2) we have that
the dimension $\mathfrak{n}(T_R) $ of the null space of $ T_R $, is given by
\[
\mathfrak{n}(T_R) = \sum_{\kappa_j \geq 1 } \kappa_j .
\]
Consider the function $ R $ multiplied by $ z^k $, which we denote by $ z^k R $.
Define for $ k = 1, 2, \ldots $ the numbers $ \mu_k $ by
\[
\mu_k = \mathfrak{n}(T_{z^{k-1} R} ) - \mathfrak{n}(T_{z^{k} R} ) .
\]
Then (see Section 2 below)
\[
\kappa_j = \# \{ k \, : \, \mu_k \geq j \}.
\]
Here $\# E $ denotes the number of elements of the set $E$.

When $R(z)$ is given in terms of a state space realization
\[
R(z)=R_0+\gamma(z I-\alpha)^{-1}\beta + zC(I-zA)^{-1}B,
\]
finding the Wiener-Hopf indices in terms of the matrices appearing in the realization
was studied in \cite{GKRb}. This approach was based on \cite{GKRa}, where the problem
was considered for the special case where $R$ takes unitary values on the unit circle. There,
significant use was made of the Douglas-Shapiro-Shields factorization of $R$, that is,
writing $R$ as $R(z)=V(z)W^*(z)$, where $V$ and $W$ are bi-inner
$m\times m$ matrix functions and $ W^*(z) $ is defined by
$W^*(z) = \bigl(W(\frac{1}{\overline{z}} ) \bigr)^* $.
(A function $F(z)$ is \emph{bi-inner} if $F(z)$ is analytic in the open unit disc
and almost everywhere unitary on the unit circle.)
The Wiener-Hopf indices were given in \cite{GKRa} in terms of realizations of $V$ and $W$, based
on earlier work in \cite{FK}.

Finding the Wiener-Hopf indices in terms of matrices in a realization is a problem that has already
some history, see \cite{BGKOT12,BGK4OT21,BGKOT21,BGKROT200,GKS,GLeR}.
Our aim in this paper is to revisit the result of \cite{GKRa} with a different method, which is more
operator theoretic, and which leads to alternative and simpler formulas for the Wiener-Hopf
indices of a rational matrix function  that  takes unitary values on the unit circle.

We conclude this introduction with a short description of the
various sections of this paper.
In Section 2, we introduce the functions and their realizations, and present the main result in Theorem
\ref{mainthm00}. Section 3 is concerned with the Toeplitz and Hankel operators corresponding to bi-inner rational matrix functions. Section 4 gives more detailed results on the unimodular function $R$, its Toeplitz operator $T_R$ and factorization. In subsection 4.1 we specify the results for the case where the bi-inner functions are scalar valued Blaschke products. We derive the main results for this special case. Section 5 is devoted to the proof of the main theorem. In Section 6 a direct connection is made with the results of \cite{GKRa}. Section 7 is an appendix with some useful observations on unitary lower triangular operators.


\section{The main result.}
To present our method to compute the Wiener-Hopf indices,
let us fix some notation.
Recall that $\{A \mbox{ on } \mathcal{X},B,C,D\}$ is a {\it realization of a function
$\Theta(z)$} if
\[
\Theta(z) = D +z C(I-zA)^{-1}B
\]
in some neighborhood of the origin.
Here $A$ is an operator on $\mathcal{X}$ and $B$ maps $\mathcal{U}$ into $\mathcal{X}$,
while $C$ maps $\mathcal{X}$ into $\mathcal{Y}$ and
$D$ maps $\mathcal{U}$ into $\mathcal{Y}$.
Two state space realizations $\{A \mbox{ on }\  \mathcal{X},B,C,D\}$
 and $\{A_1 \mbox{ on } \mathcal{X}_1 ,B_1,C_1,D_1\}$
are \emph{unitarily equivalent}  if $D = D_1$ and there
exists a unitary operator $U$ mapping $\mathcal{X}_1$ onto $\mathcal{X}$ such that
\[
A U = UA_1 \quad \mbox{and}\quad B = UB_1
\quad \mbox{and}\quad CU = C_1.
\]
Throughout we will be dealing with finite dimensional realizations, that is,
realizations of the form $\{A \mbox{ on } \mathcal{X} ,B,C,D\}$ where
the state space $\mathcal{X}$ is finite dimensional. The realization
$\{A,B,C,D\}$ is \emph{stable,} if all the eigenvalues for $A$ are contained
in the open unit disc.
Finally, we say that the realization
$\{A,B,C,D\}$ is {\it unitary } if its systems operator
\[
T = \begin{bmatrix}
      A & B \\
      C & D \\
    \end{bmatrix}:\begin{bmatrix}
      \mathcal{X} \\
      \mathcal{U}\\
      \end{bmatrix}\rightarrow
      \begin{bmatrix}
      \mathcal{X} \\
      \mathcal{Y} \\
    \end{bmatrix}
\]
is a unitary operator.

We say that $\Theta$ is a {\it bi-inner function} if
$\Theta$ is a function in $H^\infty(\mathcal{E},\mathcal{E})$ and
$\Theta(e^{i \omega})$ is almost everywhere a unitary
operator on $\mathcal{E}$.
(Here $H^\infty(\mathcal{E},\mathcal{E})$  is the Banach space
of all uniformly bounded analytic functions in the
open unit disc whose values are linear operators on $\mathcal{E}$; see \cite{sznfbk}.)
We are now ready to
present the following classical result which is
a special case of the Sz.-Nagy-Foias model theory
or characteristic function theory.

\begin{theorem}\label{thm-inner}  Let $\Theta$ be a rational function in $H^\infty(\mathcal{E},\mathcal{E})$.
Then $\Theta$ is bi-inner if and only if $\Theta$ admits
a stable unitary realization $\{A,B,C,D\}$. In this
case, all stable unitary realizations of $\Theta$ are
unitarily equivalent.
\end{theorem}

In this note, we only need the case when $\Theta$ is rational. However,
it is emphasized that Theorem \ref{thm-inner} holds in the infinite dimensional case.
The only modification is that in the infinite dimensional case,
we say that $A$ is \emph{stable} if both $A$ and $A^*$ are both pointwise stable.
For further results see the Sz.-Nagy-Foias model theory in
\cite{sznf} or its second edition  \cite{sznfbk}.

Due to the Douglas-Shapiro-Shields factorization, when computing
the\\ Wiener-Hopf indices of $R$, without loss of generality, one can assume that
$R = VW^*$ where $V$ and $W$ are two rational
bi-inner functions in $H^\infty(\mathcal{E},\mathcal{E})$.
Let
 $\{A_v \mbox{ on } \mathcal{X}_v, B_v,C_v,D_v \}$
and $\{A_w \mbox{ on } \mathcal{X}_w, B_w,C_w,  D_w\}$  be
two stable unitary realizations of $V$ and $W$ respectively.
In particular,
\begin{align}
V(z) &= D_v + z C_v (I - z A_v )^{-1} B_v, \label{defV00} \\
W(z) &= D_w + z C_w (I - z A_w )^{-1} B_w. \label{defW00}
\end{align}

Let $ -\kappa_1 , -\kappa_2  , \ldots, -\kappa_p  $ with
$ \kappa_1 \geq \kappa_2 \geq \cdots \geq \kappa_p  $ be the negative
Wiener-Hopf indices of the function $ R $
and $ R(z) = W_-(z) D(z) W_+(z) $ the Wiener-Hopf factorization of $ R $.
Then (e.g., see \cite{GGKOT63}  Theorem XXIV.4.2) we have that
the dimension $\mathfrak{n}(T_R) $ of the null space of $ T_R $, is given by
\[
\mathfrak{n}(T_R) = \sum_{\kappa_j \geq 1 } \kappa_j .
\]
Consider the function $ R $ multiplied by $ z^k $, which we denote by $ z^k R $.
Since $ z^k R(z) = W_-(z) \bigl( z^k D(z) \bigr) W_+(z)$,
the Wiener-Hopf indices of $z^k R $ are each $k$ higher  than
the corresponding  index of  $R$.
Therefore
\[
\mathfrak{n}(T_{z^k R} ) = \sum_{\kappa_j \geq k+1 }  (\kappa_j -k).
\]
Define for $ k = 1, 2, \ldots $, the numbers $ \mu_k $ by
\begin{equation}\label{defmuk10}
\mu_k = \mathfrak{n}(T_{z^{k-1} R} ) - \mathfrak{n}(T_{z^{k} R} ) =
\# \{ j \, : \, \kappa_j \geq k \}.
\end{equation}
Then (see \cite{GKvS1995} Proposition III.4.1)
\begin{equation}\label{defkappaj}
\kappa_j = \# \{ k \, : \, \mu_k \geq j \}.
\end{equation}

\medskip
The main result is the following theorem.
\medskip

\begin{theorem}\label{mainthm00}
Assume that $R = V W^*$ where
 $V$ and $W$ are two bi-inner rational functions in $H^\infty(\mathcal{E},\mathcal{E})$.
Let $\{A_v, B_v,C_v,D_v \}$
and $\{A_w, B_w,C_w,  D_w\}$  be
stable unitary realizations of $V$ and $W$,  respectively.
Let $\Omega$ be the unique solution of the Stein equation
\begin{equation}\label{lyapom00}
\Omega = A_v \Omega A_w^* + B_v B_w^*.
\end{equation}
Let $C_\circ$ be the operator mapping $\mathcal{X}_w$ into $\mathcal{E}$ defined by
\begin{equation}  
C_\circ = D_v B_w^* + C_v \Omega A_w^*.
\end{equation}
Finally,  let $Q$ be the unique solution to the Stein equation
\begin{equation}\label{defQ00}
 Q = A_w  Q A_w^* + C_\circ^* C_\circ.
\end{equation}
Then the following holds:
\begin{enumerate}
  \item The operator $Q$ is a positive contraction.
  \item The multiplicity of 1 as an eigenvalue of $ Q $  equals $\mathfrak{n}(T_R)$.
In other words, $\mathfrak{n}(T_R) = \mathfrak{n}(I-Q)$.
Moreover, for $k=0,1,2,\cdots, p$, we have
\begin{equation}\label{alphathm00}
\mathfrak{n}(T_{z^k R})  = \dim \left( \kr (I - A_w^{k} Q A_w^{*k} )\right).
\end{equation}
  \item For $k=1,2,\cdots$, consider  the sequence of integers
\begin{equation}\label{defmuk00}
\mu_k = \mathfrak{n}(I - A_w^{k-1} Q (A_w^*)^{k-1} )
- \mathfrak{n} (I - A_w^{k} Q A_w^{*k} ) .
\end{equation}
Then the negative Wiener-Hopf indices $ -\kappa_1 , \ldots , -\kappa_p $
of $ T_R $ are given by
\begin{equation}\label{defmuk001}
\kappa_j = \# \{ k : \mu_k \geq j \}, \quad (j = 1 , \ldots, p = \mu_1 ).
\end{equation}
\end{enumerate}
\end{theorem}

Notice that, once Parts 1 and 2 are proven,  Part 3 follows from
the equations \eqref{defmuk10} and \eqref{defkappaj}.

The dual statement for the positive Wiener-Hopf indices is obtained by applying
the above theorem to the function $ R^*(z)= \bigl( R(\frac{1}{\bar{z}}) \bigr)^* $; see Corollary~\ref{maincol} below.


\setcounter{equation}{0}
\section{Unitary functions on the unit circle}

Let $H^\infty(\mathcal{E},\mathcal{E})$ be the Hardy space
consisting of the set of all operator valued functions
$\Theta(z)$ on $\mathcal{E}$  analytic in the open unit disc
such that
\[
\|\Theta\|_\infty = \sup\{\|\Theta(z)\|: |z|<1\} <\infty.
\]
Let $\Theta(z) = \sum_0^\infty z^n \Theta_n$ be
the Taylor series expansion for a function in
$\Theta$ in $H^\infty(\mathcal{E},\mathcal{E})$.
Then $T_\Theta$ is the Toeplitz operator on $\ell_+^2(\mathcal{E})$
and $H_\Theta$ is the Hankel  operator on $\ell_+^2(\mathcal{E})$
defined by
\begin{equation}\label{deftoep}
T_\Theta = \begin{bmatrix}
                   \Theta_0 & 0 & 0 & \cdots \\
                   \Theta_1 & \Theta_0 & 0 & \cdots \\
                   \Theta_2 & \Theta_1 & \Theta_0 & \cdots \\
                   \vdots &  \vdots & \vdots  & \vdots \\
             \end{bmatrix}: \ell_+^2(\mathcal{E}) \rightarrow\ell_+^2(\mathcal{E}) ,
\end{equation}
\begin{equation} \label{defHank}
 H_\Theta =
\begin{bmatrix}
                \Theta_1 & \Theta_2 & \Theta_3 & \cdots \\
                \Theta_2 & \Theta_3 & \Theta_4 & \cdots \\
                 \Theta_3 & \Theta_4 & \Theta_5 & \cdots \\
                \vdots & \vdots  & \vdots & \vdots \\
              \end{bmatrix} \ell_+^2(\mathcal{E}) \rightarrow\ell_+^2(\mathcal{E}).
\end{equation}
Throughout $\widetilde{\Theta}$ is the function in
$H^\infty(\mathcal{E},\mathcal{E})$ defined by
$\widetilde{\Theta}(z) = \Theta( \bar{z} )^*$
for all $z$ with $ |z| \leq 1$. The
Taylor series expansion for $\widetilde{\Theta}(z)$ is given by
\begin{equation}\label{taylor-tilde}
\widetilde{\Theta}(z) = \sum_{n=0}^\infty z^n \Theta_n^*
\qquad (\mbox{for } |z| <1).
\end{equation}
In particular, this implies that the Hankel matrix
$H_{\widetilde{\Theta}}=H_\Theta^*$.
Finally, if $\mathcal{E} = \mathbb{C}$, then
$H^\infty(\mathbb{C},\mathbb{C})$ is simply denoted by $H^\infty$.

As before, let $\Theta$ be a bi-inner function in
$H^\infty(\mathcal{E},\mathcal{E})$.
In this case, the corresponding Toeplitz matrix $T_\Theta$ is an isometry
on $\ell_+^2(\mathcal{E})$.
Because $\widetilde{\Theta}$ is also bi-inner,
$T_{\widetilde{\Theta}}$ is also an isometry on $\ell_+^2(\mathcal{E})$.
Let  $\mathfrak{H}(\Theta)$ and $\mathfrak{H}(\widetilde{\Theta})$
denote the orthogonal complements of the ranges of $T_\Theta$ and
$ T_{\widetilde{\Theta}} $, respectively, that is,
\begin{equation}\label{defHtheta}
\mathfrak{H}(\Theta)=\ell_+^2(\mathcal{E})\ominus T_\Theta \ell_+^2(\mathcal{E})
\quad \mbox{and}\quad
 \mathfrak{H}(\widetilde{\Theta}) = \ell_+^2(\mathcal{E})\ominus T_{\widetilde{\Theta}} \ell_+^2(\mathcal{E}).
\end{equation}
By consulting the Appendix, the Hankel operator
$H_\Theta$ can be viewed as a unitary operator
from $\mathfrak{H}(\widetilde{\Theta})$ onto $\mathfrak{H}(\Theta)$.
In particular, $\mathfrak{H}(\Theta)$ equals the range of
$H_\Theta$ and $\mathfrak{H}(\widetilde{\Theta})$ equals the range of
$H_\Theta^* = H_{\widetilde{\Theta}}$. Finally, the corresponding orthogonal projections
are given by
\begin{equation}\label{PH0}
P_{_{\mathfrak{H}( \Theta)}} = H_\Theta H_\Theta^* =
I - T_\Theta T_\Theta^*
\quad \mbox{and}\quad
P_{_{\mathfrak{H}(\widetilde{\Theta})}} = H_\Theta^* H_\Theta =
I - T_{\widetilde{\Theta}} T_{\widetilde{\Theta}}^*.
\end{equation}

Let $\{A \mbox{ on } \mathcal{X},B,C,D\}$ be any stable unitary
realization of a rational  bi-inner function $\Theta$ in $H^\infty(\mathcal{E},\mathcal{E})$.
Then its \emph{observability operator} $\Gamma $ and 
\emph{controllability operator} $\Upsilon$ are the 
 defined by
\begin{equation}\label{conaobs}
\Gamma = \begin{bmatrix}
                    C \\
                    CA \\
                    C A^2 \\
                    \vdots \\
                  \end{bmatrix}:\mathcal{X} \rightarrow\ell_+^2(\mathcal{E})
  \mbox{ and }
\Upsilon = \begin{bmatrix}
                    B &
                    AB &
                    A^2 B &
                    \cdots &
                  \end{bmatrix}: \ell_+^2(\mathcal{E})\rightarrow\mathcal{X}.
\end{equation}
Because $\begin{bmatrix}
           C^* & A^* \\
         \end{bmatrix}^*$ is an isometry and $A$ is  stable,
the observability operator $\Gamma$ is an isometry.
Likewise, since $\begin{bmatrix}
           B & A \\
         \end{bmatrix}^*$ is an isometry and $A^*$ is  stable,
the controllability  operator $\Upsilon $ is a co-isometry.
Using the fact that $\Theta_n = C A^{n-1}B$ for all integers
$n \geq 1$, it follows that the Hankel operator $H_\Theta$ admits a
factorization of the form:
\begin{equation}\label{hankobs}
 H_\Theta = \Gamma \Upsilon.
\end{equation}
Above we noticed that $H_\Theta$ can be viewed as a unitary
operator from $\mathfrak{H}(\widetilde{\Theta})$
onto $\mathfrak{H}(\Theta)$.
Since $\Gamma $ is an isometry and
$\Upsilon $ is a co-isometry, the equalities
$\mathfrak{H}(\Theta) = \im (\Gamma )$
and $\mathfrak{H}(\widetilde{\Theta}) = \im (\Upsilon^*)$ hold.
(Here $\im$ denotes the range of an operator.)
The equation $H_\Theta = \Gamma \Upsilon $ with \eqref{PH0},   readily
implies that
\begin{equation}\label{GammaH}
P_{_{\mathfrak{H}(\Theta)}} = \Gamma \Gamma^*
\quad \mbox{and}\quad
P_{_{\mathfrak{H}(\widetilde{\Theta})}} = \Upsilon^* \Upsilon.
\end{equation}


\setcounter{equation}{0}
\section{The function $R= VW^*$.  }
Let $V$ and $W$ be two rational bi-inner functions in
$H^\infty(\mathcal{E},\mathcal{E})$. Let $R$ be the rigid function in
$L^\infty(\mathcal{E},\mathcal{E})$ defined by
\[
R(e^{i\omega}) = V(e^{i\omega}) W(e^{i\omega})^* \qquad (\mbox{for } 0 \leq \omega \leq 2\pi).
\]
(Recall that a function $\Xi$ in $L^\infty(\mathcal{E},\mathcal{E})$ is rigid if
$\Xi(e^{i \omega})$ is almost everywhere a unitary operator on $\mathcal{\mathcal{E}}$.)
Let $T_R$ be the Toeplitz operator on $\ell_+^2(\mathcal{E})$ determined by $R$, that is,
if $R(e^{i \omega}) = \sum_{-\infty}^\infty e^{i n \omega} R_n$ is the Fourier series expansion
for $R$, then
\[
T_R = \begin{bmatrix}
        R_0 & R_{-1} & R_{-2} & \cdots \\
        R_1 & R_{0} & R_{-1} & \cdots \\
        R_2 & R_{1} & R_{0} & \cdots \\
        \vdots  & \vdots & \vdots & \ddots \\
      \end{bmatrix}:\ell_+^2(\mathcal{E}) \rightarrow \ell_+^2(\mathcal{E}).
\]
Since $R = VW^*$, it follows that
\begin{equation}\label{toephank}
 T_R = T_V T_W^* + H_V H_W^*.
\end{equation}
Because $V$ and $W$ are bi-inner,
 $H_V$ is a unitary operator from
$\mathfrak{H}(\widetilde{V})$ onto $\mathfrak{H}(V)$,
and $H_W$ is a unitary operator from
$\mathfrak{H}(\widetilde{W})$ onto $\mathfrak{H}(W)$.
This readily implies that
\begin{equation}\label{toephankY}
 T_R = T_V T_W^* + H_V Y H_W^*
\end{equation}
where $Y$ is the contraction mapping
$\mathfrak{H}(\widetilde{W})$ into
$\mathfrak{H}(\widetilde{V})$ defined by
\begin{equation}\label{defY}
Y = P_{_{\mathfrak{H}(\widetilde{V})}}|\mathfrak{H}(\widetilde{W})
: \mathfrak{H}(\widetilde{W}) \rightarrow \mathfrak{H}(\widetilde{V}).
\end{equation}
Since $V$ and $W$ are both bi-inner,
\[
\ell_+^2(\mathcal{E}) = \im (T_V) \oplus \im (H_V)
\quad \mbox{and}\quad
\ell_+^2(\mathcal{E}) = \im (T_W) \oplus \im (H_W).
\]
Using this with  $T_R = T_VT_W^* + H_VYH_W^*$, we see that
$T_R$ admits a "singular value type" decomposition of the form:
\begin{equation}\label{svdTR}
T_R = T_VT_W^* + H_VYH_W^*
= \begin{bmatrix}
    T_V & H_V \\
  \end{bmatrix}\begin{bmatrix}
                 I & 0 \\
                 0 & Y \\
               \end{bmatrix}\begin{bmatrix}
    T_W^* \\ H_W^* \\
  \end{bmatrix}.
\end{equation}
Here
\[
\begin{bmatrix}
    T_V & H_V \\
  \end{bmatrix}:\begin{bmatrix}
    \ell_+^2(\mathcal{E}) \\ \mathfrak{H}(\widetilde{V}) \\
  \end{bmatrix}\rightarrow\ell_+^2(\mathcal{E})
  \quad \mbox{and}\quad
  \begin{bmatrix}
    T_W^* \\ H_W^* \\
  \end{bmatrix}:\ell_+^2(\mathcal{E})
  \rightarrow \begin{bmatrix}
    \ell_+^2(\mathcal{E}) \\ \mathfrak{H}(\widetilde{W}) \\
  \end{bmatrix}
\]
are both unitary operators. Moreover, the
middle term
\begin{equation}\label{middle}
\begin{bmatrix}
                 I & 0 \\
                 0 & Y \\
               \end{bmatrix}  = \begin{bmatrix}
                 I & 0 \\
                 0 & P_{_{\mathfrak{H}(\widetilde{V})}}|\mathfrak{H}(\widetilde{W}) \\
               \end{bmatrix}:
\begin{bmatrix}
    \ell_+^2(\mathcal{E}) \\ \mathfrak{H}(\widetilde{W}) \\
  \end{bmatrix}
\rightarrow \begin{bmatrix}
    \ell_+^2(\mathcal{E}) \\ \mathfrak{H}(\widetilde{V}) \\
  \end{bmatrix}
\end{equation}
is a contraction.

Due to the decomposition of $T_R$ in \eqref{svdTR}, it follows that all the
properties such as invertibility and Fredholmness
of the operator $T_R$ are the same as those of the contraction $Y$.

It is noted that $x$ is in $\kr(Y)$ if and only if
$x$ is in $\mathfrak{H}(\widetilde{W})$ and $P_{_{\mathfrak{H}(\widetilde{V})}} x = 0$,
or equivalently, $x$ is in $\mathfrak{H}(\widetilde{W})$ and
$x$ is in $\im (T_{\widetilde{V}}) = \mathfrak{H}(\widetilde{V})^\perp$.
In other words,
\[
\kr(Y) = \im (T_{\widetilde{V}})\bigcap \mathfrak{H}(\widetilde{W})
\quad \mbox{and} \quad
\kr(Y^*) = \im (T_{\widetilde{W}})\bigcap \mathfrak{H}(\widetilde{V}),
\]
where the second equality follows from a similar argument.

Recall that an  operator $T$ mapping $\mathcal{X}$ into $\mathcal{Y}$
admits a Moore-Penrose inverse $T^{pinv}$  if the operator $T\vert\kr(T)^\perp$ mapping
$\kr(T)^\perp$ into the range of $T$ is invertible. In this case, the
Moore-Penrose inverse of $T$ is given by $T^{pinv}=\left(T\vert\kr(T)^\perp\right)^{-1}P_{\im (T)}$.
By consulting the form of $T_R$ in \eqref{svdTR},
we obtain the following result.

\begin{proposition}\label{prop-R}
Let $R=VW^*$ where $V$ and $W$ are both bi-inner functions in
$H^\infty(\mathcal{E},\mathcal{E})$. Moreover, let $Y$
be the contraction mapping $\mathfrak{H}(\widetilde{W})$ into
$\mathfrak{H}(\widetilde{V})$ defined by
$Y = P_{_{\mathfrak{H}(\widetilde{V})}}|\mathfrak{H}(\widetilde{W})$.
Then the following holds.
\begin{enumerate}
  \item The operator $T_R$ is invertible if and only if
  $Y$ is invertible. In this case,
  \begin{equation}\label{invTR}
    T_R^{-1} = T_W T_V^* + H_W Y^{-1}H_V^*.
  \end{equation}
  \item The subspaces  $\kr(T_R)$ and $\kr(Y)$ have the same dimension.
  In fact,
  \begin{equation}\label{kerTR}
    \kr(T_R) = H_W \kr(Y) \quad \mbox{and}\quad
    \kr(Y) = \im (T_{\widetilde{V}})\bigcap \mathfrak{H}(\widetilde{W}).
  \end{equation}
  \item The subspaces  $\kr(T_R^*)$ and $\kr(Y^*)$ have the same dimension.
  In particular,
  \begin{equation}\label{kerTR1}
    \kr(T_R^*) = H_V \kr(Y^*)
    \quad \mbox{and}\quad
    \kr(Y^*) = \im (T_{\widetilde{W}})\bigcap \mathfrak{H}(\widetilde{V}).
  \end{equation}
  \item The subspaces $\im (T_R)^\perp$ and $\im (Y)^\perp$
  have the same dimension. In fact,
  \begin{equation}\label{kerTR2}
   \im (T_R)^\perp =  \kr(T_R^*) = H_V \kr(Y^*) = H_V \,\im (Y)^\perp.
  \end{equation}
  \item The operator $T_R$ admits a Moore-Penrose restricted inverse if and only if
  $Y$ admits a Moore-Penrose restricted inverse. In this case,
  \begin{equation}\label{MPenTR}
    T_R^{pinv} = T_W T_V^* + H_W Y^{pinv}H_V^*.
  \end{equation}
\end{enumerate}
\end{proposition}


\subsection{The Blaschke product case}

In this section, to gain some insight into the general case,  we will study the
contraction $Y=P_{_{\mathfrak{H}(\varphi)}}|\mathfrak{H}(m)$
mapping $\mathfrak{H}(m)$ into $\mathfrak{H}(\varphi)$ when
$m$ and $\varphi$ are Blaschke  products.
We say that a function $b(z)$ is a \emph{Blaschke product}  if
\begin{equation}\label{blaschke}
b(z) = \zeta \prod_{k=1}^n \frac{z-\alpha_k}{1-\overline{\alpha}_k z}
\qquad (\mbox{where } |\alpha_k| < 1 \mbox{ for all } k).
\end{equation}
(Here $\zeta$ is a complex number on the unit circle.)
Moreover, $n = \deg(b)$ is the \emph{degree}  of the Blaschke product.
Throughout we will only consider Blaschke products of finite degree.
So if we say that $b(z)$ is a Blaschke product,
then we assume that  $b(z)$  is a function of the form \eqref{blaschke}
and the degree of $b$ is finite. It is well known that
$b(z)$ is a rational inner function in $H^\infty$ if and only if
$b$ is a Blaschke product (of finite degree).
Furthermore,   $b(z)$ is a Blaschke product of
degree $n$ if and only if $b(z)$ admits a stable unitary
realization $\{A \mbox{ on } \mathcal{X}, B ,C ,D\}$ where
$n$ is the dimension of the state space $\mathcal{X}$.
Moreover, in this case, the zeros of $b(z)$ are precisely the
eigenvalues of $A$. These are classical results in
systems theory and are a special case of the
Sz.-Nagy-Foias model theory for $C_0$ contractions; see  \cite{sznf,sznfbk}.

The Fourier transform $\mathfrak{F}$ is the unitary operator mapping
$\ell_+^2$ onto $H^2$ defined by
\begin{equation}\label{fourier}
\mathfrak{F}\begin{bmatrix}
              x_0 & x_1  & x_2 & \cdots \\
            \end{bmatrix}^{tr} = \sum_{k=0}^\infty x_n z^n
            \qquad (\begin{bmatrix}
              x_0 & x_1  & x_2 & \cdots \\
            \end{bmatrix}^{tr} \in \ell_+^2).
\end{equation}
(The transpose of a vector  is denoted by $^{tr}$.)
Assume that $\{A \mbox{ on } \mathcal{X}, B ,C ,D\}$  is a
stable unitary realization for a Blaschke product $b(z)$
of degree $n$.
Let $\Gamma$ mapping $\mathcal{X}$ into $\ell_+^2$ be the
observability operator formed by the pair $\{C,A\}$.
Recall that $\Gamma$ is an isometry. Moreover, the
range of $\Gamma$ equals $\mathfrak{H}(b)$.
Because the degree of the state space is $n$, it follows
that the dimension of $\mathfrak{H}(b)$
(denoted by $\dim(\mathfrak{H}(b))$) equals $n$.
Moreover, the Fourier transform of $\mathfrak{H}(b)$ is given by
\begin{equation}\label{FHb}
\mathfrak{F}\big(\mathfrak{H}(b)\big) = \{\mathfrak{F}\big( \Gamma x\big) :x\in \mathcal{X}\}
= \{C(I - z A)^{-1} x : x \in \mathcal{X}\}.
\end{equation}
Because the pair $\{C,A\}$ is observable,
\begin{equation}\label{FHb1}
\mathfrak{F}\big(\mathfrak{H}(b)\big) = \left\{\frac{p(z)}{\det[I-zA]}: p(z) \mbox{ is a polynomial
of degree} < \dim({\mathcal{X}})\right\}.
\end{equation}
So if $b(z)$ is the Blaschke product of degree $n$ given in
\eqref{blaschke}, then
\begin{equation}\label{FHb2}
\mathfrak{F}\big(\mathfrak{H}(b)\big) = \left\{\frac{p(z)}{\prod_{k=1}^n(1-\overline{\alpha}_k z)}: p(z) \mbox{ is a polynomial
of degree} < n\right\}.
\end{equation}

If $A$ is a contraction on $\mathcal{X}$, then $D_A = (I - A^*A)^{\frac{1}{2}}$ is the
\emph{ defect operator} for $A$. Moreover, $\mathfrak{D}_A$ is the
range of $D_A$. Finally,
$\mathfrak{d}_A = \dim(\mathfrak{D}_A)$ is the \emph{ defect index} for $A$.

If $A$ is a stable contraction on a finite dimensional space $\mathcal{X}$
and $\psi$ is a function in $H^\infty$, then
 $\psi(A)$ is the function of $A$ defined by
 \[
 \psi(A) = \sum_{n=0}^\infty \psi_n A^n
 \qquad \mbox{where}\qquad
 \psi(z) = \sum_{n=0}^\infty \psi_n z^n.
 \]
In this case,   $\|\psi(A)\| \leq \|\psi\|_\infty$.
For further results in the infinite dimensional setting, see the function theory in \cite{sznf,sznfbk}.
 This sets the stage for the following result.

\begin{lemma}\label{lem-mphi} Let $m(z)$  and $\varphi(z)$ be 
finite  Blaschke products
in $H^\infty$.
Consider the contraction $Y=P_{_{\mathfrak{H}(\varphi)}}|\mathfrak{H}(m)$
mapping $\mathfrak{H}(m)$ into $\mathfrak{H}(\varphi)$.
Let $\{A,B,C,D\}$ be a unitary stable realization for $m$.
Then the following holds.
\begin{enumerate}
  \item There exists a unitary operator $\Psi$ mapping the range of
  $Y$ onto $\mathfrak{D}_{_{\widetilde{\varphi}(A)}}$ such that
   \begin{equation}\label{defectPsi}
     \Psi P_{_{\mathfrak{H}(\varphi)}}\Gamma_m = D_{_{\widetilde{\varphi}(A)}}.
   \end{equation}
   (Here $\Gamma_m$ is the observability operator formed by the state space
   realization $\{A \mbox{ on } \mathcal{X},B,C,D\}$ for $m$.)
   In this case,
   \begin{equation}\label{defect00}
    \mathfrak{d}_{\widetilde{\varphi}(A)} =  \mathfrak{d}_{ \varphi (A)} =\min\left\{\deg(\varphi),\deg(m)\right\}  =
    \rank\big( P_{_{\mathfrak{H}(\varphi)}}\vert \mathfrak{H}(m)\big).
    \end{equation}
  \item In particular,  $\deg(\varphi) \leq \deg(m)$ if and only if  the range of
   the contraction $Y=P_{_{\mathfrak{H}(\varphi)}}|\mathfrak{H}(m)$ equals
   $\mathfrak{H}(\varphi)$. In this case,
   $\mathfrak{d}_{\varphi(A)} =
   \mathfrak{d}_{\widetilde{\varphi}(A)} =\deg(\varphi)$, and
   \begin{equation}\label{blaschke00}
   \dim(\kr(Y)) = \deg(m)-\deg(\varphi) \qquad (\mbox{when } \deg(\varphi) \leq \deg(m)).
   \end{equation}
  \item The operator $Y$ is one to one if and only if
  $\deg(m) \leq \deg(\varphi)$. In this case,
  $\mathfrak{d}_{\varphi(A)} =
   \mathfrak{d}_{\widetilde{\varphi}(A)} =\deg(m)$,
   \begin{equation}\label{dimmphi}
    \dim(\kr(Y^*)) = \dim(\im (Y)^\perp) = \deg(\varphi) - \deg(m).
\end{equation}
  \item If $\deg(\varphi) < \deg(m)$, then the Blaschke product
  \begin{equation}\label{aak}
   \varphi(z) = \frac{C(I-z A)^{-1}x}{C(I-z A)^{-1}\widetilde{\varphi}(A)x}
    \qquad (\mbox{if } 0 \neq x \in \mathfrak{D}_{_{\widetilde{\varphi}(A)}}^\perp).
  \end{equation}
\end{enumerate}
\end{lemma}

\noindent {\sc Proof.}
Because  $\{A \mbox{ on } \mathcal{X},B,C,D\}$ is a
stable unitary realization of $m$, we have $\deg(m) = \dim(\mathcal{X})$.
Recall that the Hankel operator
$H_m = \Gamma_m \Upsilon_m$ where $\Gamma_{m}$
mapping $\mathcal{X}$ into $\ell_+^2$ is the observability
operator formed by $\{C,A\}$. Moreover, $\Upsilon_m $ mapping $\ell_+^2$
onto $\mathcal{X}$
is the controllability operator determined  by $\{A,B\}$. Furthermore,
$\Gamma_{m}$ is an isometry and $\Upsilon_m $ is a co-isometry.
Since the range of $H_m$ equals $\mathfrak{H}(m)$, it follows
that the subspace $\mathfrak{H}(m)$ equals the range of $\Gamma_{m}$.
Notice that $h$ is in $\mathfrak{H}(m)$ if and only if
$h = \Gamma_{m} x$ for some $x$ in $\mathcal{X}$.
In fact, this $x$ is uniquely determined by $h$ and given by
$x = \Gamma_{m}^* h$.
Using  $P_{_{\mathfrak{H}(\varphi)}}  = I - T_\varphi T_\varphi^*$
with $T_\varphi^* \Gamma_{m} = \Gamma_{m} \widetilde{\varphi}(A)$
and $ T_\varphi $ is an isometry, we have
\begin{align}\label{long}
\|P_{_{\mathfrak{H}(\varphi)}} \Gamma_{m} x\|^2 &=
\|(I - T_\varphi T_\varphi^*)\Gamma_{m} x\|^2 =
\|\Gamma_{m} x\|^2 - \|T_\varphi T_\varphi^*\Gamma_{m} x\|^2\nonumber\\
&= \|x\|^2 - \| T_\varphi^*\Gamma_{m} x\|^2 =
\|x\|^2 - \|\Gamma_{m} \widetilde{\varphi}(A) x\|^2\nonumber\\
&=\|x\|^2 - \|\widetilde{\varphi}(A) x\|^2 =
 \langle x,(I - \widetilde{\varphi}(A)^*\widetilde{\varphi}(A))x \rangle   \nonumber\\
&= \|(I - \widetilde{\varphi}(A)^*\widetilde{\varphi}(A))^{\frac{1}{2}}x\|^2
=\|D_{_{\widetilde{\varphi}(A)}}x\|^2.
\end{align}
Hence  there
exists a unitary operator $\Psi$ mapping the range of
  $Y$ onto $\mathfrak{D}_{_{\widetilde{\varphi}(A)}}$ such that
   \begin{equation}\label{defectPsi54}
     \Psi P_{_{\mathfrak{H}(\varphi)}}\Gamma_{m} = D_{_{\widetilde{\varphi}(A)}}.
 \end{equation}
This proves equation \eqref{defectPsi} in Part 1.

Now let us show  that $\mathfrak{d}_{_{\widetilde{\varphi}(A)}} = \min\{\deg(\varphi),\deg(m)\}$.
To this end, first assume that $\deg(\varphi) \leq \deg(m)$.
Then we claim that $Y$ is onto $\mathfrak{H}(\varphi)$, and thus,
the rank of $Y$ equals  $\dim(\mathfrak{H}(\varphi)) = \deg(\varphi)$.
Assume that a vector $h\in \mathfrak{H}(\varphi)$ is orthogonal to the range of
$Y=P_{_{\mathfrak{H}(\varphi)}}|\mathfrak{H}(m)$.
Then it follows from
\eqref{kerTR1} with $ \widetilde{V} = \varphi $ and $ \widetilde{W}=m $ that
$h$ is also a vector in the range of $T_m$, that is,
$h\in \mathfrak{H}(\varphi)\cap T_m \ell_+^2 $.
By consulting \eqref{FHb1} or \eqref{FHb2}, we see that the Fourier transform
$\mathfrak{F}\big(\mathfrak{H}(\varphi)\big)$ of the subspace
$\mathfrak{H}(\varphi)$,  consists of a set of rational functions,
with at most $\deg(\varphi)-1$ zeros. The Fourier transform of $h$ is given by
\[
\widehat{h}(z) =  \big(\mathfrak{F} h\big)(z) \in
\Big(\mathfrak{F}\big(\mathfrak{H}(\varphi)\big)\cap m H^2\Big).
\]
Since $\deg(\varphi) \leq \deg(m)$, and $m$ is rational with
$\deg(m)$ zeros, the subspace
$\Big(\mathfrak{F}\big(\mathfrak{H}(\varphi)\big)\cap m H^2\Big) = \{0\}$.
Therefore $h=0$ and the operator $Y$ is onto, whenever
$\deg(\varphi)\leq\deg(m)$.
This with
$\Psi P_{_{\mathfrak{H}(\varphi)}}\Gamma_{m} = D_{_{\widetilde{\varphi}(A)}}$, implies that
$\mathfrak{d}_{_{\widetilde{\varphi}(A)}}= \deg(\varphi)$.
Replacing $\varphi$ with $\widetilde{\varphi}$ shows that
$\mathfrak{d}_{_{ \varphi (A)}}= \deg(\varphi)$ when $\deg(\varphi) \leq \deg(m)$.

Now assume that $\deg(m) \leq \deg(\varphi)$.
 Clearly,  $Y$ and $Y^*$ have the same rank.
Notice that $Y^*$ is the contraction determined by
\[Y^*  = P_{_{\mathfrak{H}(m)}}|\mathfrak{H}(\varphi):
\mathfrak{H}(\varphi) \rightarrow \mathfrak{H}(m).
\]
So $Y^*$ is the same form as $Y$, except $m$ and $\varphi$
interchange places. By our previous analysis
 $\rank(Y^*) = \deg(m)$ and $Y^*$ is onto $\mathfrak{H}(m)$.  So $Y$ is one to one.
Recall that  $\Psi P_{_{\mathfrak{H}(\varphi)}} \Gamma_{m} = D_{_{\widetilde{\varphi}(A)}}$.
Because  $Y$ is one to one, $D_{_{\widetilde{\varphi}(A)}}$ must also be one to one.
Since $D_{_{\widetilde{\varphi}(A)}}$ is one to one and $\dim(\mathcal{X}) = \deg(m)$,
we see that $\mathfrak{d}_{_{\widetilde{\varphi}(A)}} = \deg(m)$.
This completes the proof of Part 1.

To prove Part 2, we showed that
if $ \deg (\varphi) \leq \deg(m)$  then $ Y$ is onto $ \mathfrak{H} (\varphi ) $.
Moreover, $ \dim(\kr(Y)) + \dim(\im(Y)) = \dim\left( \mathfrak{H}(m) \right)$,
which proves \eqref{dimmphi}.
On the other hand if $\deg(\varphi) > \deg(m)$, then $\rank Y = \rank Y^* =\deg{m} $
and hence $Y$ is not onto $ \mathfrak{H} (\varphi ) $.

Part 3 is proven in the same way by replacing $Y $ by $Y^* $.

To establish Part 4, assume that  $\deg(\varphi) < \deg(m)$.
Then there exists a nonzero $x$ such that $D_{_{\widetilde{\varphi}(A)}}x =0$, or
equivalently, $x = \widetilde{\varphi}(A)^* \widetilde{\varphi}(A)x$.
Using
\[
\Psi P_{_{\mathfrak{H}(\varphi)}}\Gamma_{_{m}}x = D_{_{\widetilde{\varphi}(A)}}x=0,
\]
we have
$P_{_{\mathfrak{H}(\varphi)}} \Gamma_{m}x =0$.
By employing
$P_{_{\mathfrak{H}(\varphi)}} = I-T_\varphi T_\varphi^*$,  we obtain
\[
0 = P_{_{\mathfrak{H}(\varphi)}} \Gamma_m x
= (I-T_\varphi T_\varphi^*) \Gamma_m x =
\Gamma_m  x - T_\varphi  \Gamma_m  \widetilde{\varphi}(A) x.
\]
In other words,
$\Gamma_m x = T_\varphi  \Gamma_m  \widetilde{\varphi}(A) x$.
By taking the Fourier transform, we arrive at
\begin{equation}\label{scalar}
C(I  - z A)^{-1}x = \varphi(z)C(I- zA)^{-1}\widetilde{\varphi}(A) x.
\end{equation}
Since $\Gamma_{m} $ is one to one and $\widetilde{\varphi}(A)x$ is nonzero,
$C(I  - z A)^{-1}\widetilde{\varphi}(A) x$ is a nonzero  function in  $H^2$.
Hence the rational function  $C(I  - z A)^{-1}\widetilde{\varphi}(A)x$
is nonzero on the unit circle.
Dividing by $C(I- zA)^{-1}\widetilde{\varphi}(A) x$, yields the
formula that we have been looking for, that is,
\[
\varphi(z) = \frac{C(I- zA)^{-1} x}{C(I- zA)^{-1}\widetilde{\varphi}(A) x}.
\]
This completes the proof.
\epr

\begin{remark} Let $R$ be the rational rigid function in $L^\infty$ given by
$R(z) = \varphi(z) \overline{m(1/\bar{z})}$,
where $\varphi$ and $m$ are two Blaschke products.
Let $T_R$ be the Toeplitz matrix on $\ell_+^2$ determined by $R$.
By consulting Proposition \ref{prop-R} and Lemma \ref{lem-mphi}, we readily obtain
the following.
\begin{enumerate}
  \item The operator $T_R$ is invertible if and only if
 $\deg(\varphi) = \deg(m)$.
  \item The kernel of  $T_R$ is nonzero  if and only if
 $\deg(\varphi) < \deg(m)$. In this case, the operator
$T_R$ is onto $\ell_+^2$ and $\dim(\kr(T_R)) = \deg(m) - \deg(\varphi)$.
  \item The subspace $\im(T_R)^\perp$ is nonzero  if and only if
 $\deg(m) < \deg(\varphi)$. In this case, the range of the operator
$T_R$ is closed, $\kr(T_R) = \{0\}$  and  $\dim(\im(T_R)^\perp) = \deg(\varphi) - \deg(m)$.
\end{enumerate}
\end{remark}

 \begin{corollary}\label{cor-phi}
 Let $A$ be a stable contraction on a finite dimensional space $\mathcal{X}$ with defect index one,
 and $\varphi$ a rational Blaschke product in $H^\infty$.
 Then the defect index of $A^*$ equals one and
 \begin{equation}\label{min00}
\mathfrak{d}_{_{\varphi(A)}} = \mathfrak{d}_{_{\varphi(A^*)}}
= \min\{\deg(\varphi),\dim(\mathcal{X})\}.
\end{equation}
 \end{corollary}

 \noindent {\sc Proof.} If $A$ acts on a finite dimensional
 space, then it is clear that $A^*A$ and $AA^*$ have the same
 number of eigenvalues equal to one. Hence $A$ and $A^*$ have the
 same defect index.
 In particular, this holds for
 defect index one.
(If $\mathcal{X}$ is infinite dimensional, then
 the Sz.-Nagy-Foias characteristic function theory, shows
 that for a stable contraction $A$, the defect index
 $\mathfrak{d}_A = \mathfrak{d}_{A^*}$. However, this fact is not needed here.)

Let   $C$ be any
 operator mapping $\mathcal{X}$ into $\mathbb{C}$ such that
 $I-A^*A = C^*C$. Notice that $\begin{bmatrix}
      A^* & C^* \\
    \end{bmatrix}^*$ is an isometry from $\mathcal{X}$ into $\mathcal{X}\oplus \mathbb{C}$.
Hence there exists an     operator $B$ mapping $\mathbb{C}$ into $\mathcal{X}$ and
 a  scalar $D$  such that
\[
T = \begin{bmatrix}
      A & B \\
      C & D \\
    \end{bmatrix}:\begin{bmatrix}
      \mathcal{X} \\
      \mathbb{C}\\
      \end{bmatrix}\rightarrow
      \begin{bmatrix}
      \mathcal{X} \\
      \mathbb{C}\\
    \end{bmatrix}
\]
is a unitary operator. Therefore
$\{A,B,C,D\}$ is a  stable unitary realization of
a Blaschke product $m(z)$ with degree $\dim(\mathcal{X})$; see
Theorem \ref{thm-inner}.
 Applying Lemma \ref{lem-mphi} yields \eqref{min00}.
  \epr

  \begin{proposition}\label{prop-clt00}
  Let $A$ be a stable contraction on $\mathcal{X}$ with defect index one.
  Let $\varphi$ be a Blaschke product such that
  $\deg(\varphi) < \dim(\mathcal{X})$. Then the
  inner function $\varphi$ is given by
  \begin{equation}\label{clt}
\varphi(z) = \frac{C(I-z A)^{-1}\varphi(A^*)x}{C(I-z A)^{-1}x}
\qquad (\mbox{where } x=\varphi(A^*)^*\varphi(A^*)x \mbox{ and } x \neq 0).
\end{equation}
Here $C$ is any operator mapping $\mathcal{X}$ into $\mathbb{C}$ such that $I-A^*A = C^*C$.
If $\theta$ is any function in $H^\infty$ such that
$\theta(A) = \varphi(A)$ and $\|\theta\|_\infty \leq 1$, then $\theta(z) = \varphi(z)$.
  \end{proposition}

\noindent {\sc Proof.} Let $C$ be an operator from $\mathcal{X}$
into $\mathbb{C}$ such that $I - A^*A =C^*C$. Let $\Gamma$ be the
observability operator generated by $\{C,A\}$. Recall that
$\Gamma$ is an isometry and
$T_\varphi^* \Gamma = \Gamma \widetilde{\varphi}(A)$.
By taking the adjoint of
 $T_\varphi^* \Gamma = \Gamma \widetilde{\varphi}(A)$,
 we obtain $\Gamma^* T_\varphi =   \varphi(A^*)\Gamma^*$.
 Multiplying by $\Gamma$ on the right and left yields
$\Gamma\Gamma^* T_\varphi \Gamma=   \Gamma\varphi(A^*)$.
Hence
\begin{equation}\label{clt2}
\Gamma\Gamma^* T_\varphi \Gamma x = \Gamma\varphi(A^*) x.
\end{equation}
Here $x$ is any nonzero vector satisfying
$x=\varphi(A^*)^*\varphi(A^*)x$. In other words, $x$ is an eigenvector with eigenvalue $1$
for the operator $\varphi(A^*)^*\varphi(A^*)$.  Since $\deg(\varphi) < \dim(\mathcal{X})$,
the defect index $\mathfrak{d}_{\varphi(A^*)} = \deg(\varphi)$ and $1$ is indeed
an eigenvalue for $\varphi(A^*)^*\varphi(A^*)$; see Corollary \ref{cor-phi}.

Because $\Gamma$ is an isometry,
 $\Gamma\Gamma^*$ is an orthogonal projection.
 Using this with the fact that $T_\varphi$ is an isometry, we have
\[
\|x\| =\|T_\varphi \Gamma x\| \geq
\|\Gamma\Gamma^* T_\varphi \Gamma x\|=\|\Gamma\varphi(A^*) x\| = \|\varphi(A^*) x \| =\|x\|.
\]
(Recall that $x= \varphi(A^*)^*\varphi(A^*)x$.)
Therefore we have equality, and thus,
\[
T_\varphi \Gamma x = \Gamma\Gamma^* T_\varphi \Gamma x = \Gamma\varphi(A^*) x
\qquad (\mbox{when } 0\neq x \in \mathfrak{D}_{\varphi(A^*)}^\perp).
\]
By taking the Fourier transform, we have
\[
 \varphi(z) C(I-z A)^{-1}x  =
\mathfrak{F}\big(T_\varphi \Gamma x\big)(z) =
\mathfrak{F}\big( \Gamma \varphi(A^*) x\big)(z)
=   C(I-z A)^{-1}\varphi(A^*)x.
\]
In other words,
\[
\varphi(z) C(I-z A)^{-1}x =  C(I-z A)^{-1}\varphi(A^*)x.
\]
Dividing by the function $C(I-z A)^{-1}x$, yields
the equation that we have been looking for, that is,
\[
\varphi(z) = \frac{C(I-z A)^{-1}\varphi(A^*)x}{C(I-z A)^{-1}x}.
\]

Let us show that if  $\theta$ is a function
$H^\infty$ such that $\theta(A^*) = \varphi(A^*)$ and the $H^\infty$ norm  $\|\theta\|_\infty \leq 1$,
then $\theta(z) = \varphi(z)$.

Assume that $\theta$ is a function
$H^\infty$ such that $\widetilde{\theta}(A) = \widetilde{\varphi}(A)$
 and $\|\theta\|_\infty \leq 1$.
By taking the adjoint, we see that $\theta(A^*) = \varphi(A^*)$.
In this case,
\[
T_\theta^* \Gamma = \Gamma \widetilde{\theta}(A) =
 \Gamma \widetilde{\varphi}(A).
\]
Since $T_\varphi^* \Gamma = \Gamma \widetilde{\varphi}(A)$, we also have
$T_\varphi^* \Gamma = T_\theta^* \Gamma$, or equivalently,
$ \Gamma^* T_\varphi = \Gamma^*T_\theta$. Multiplying by
$\Gamma$ on both sides, we obtain with $x$ as above
\[
\Gamma\Gamma^*T_\theta \Gamma x =  \Gamma\Gamma^*T_\varphi \Gamma x
= \Gamma   \varphi(A^*) x.
\]
Because $\Gamma\Gamma^*$ is an orthogonal projection,
$\| \varphi(A^*) x\|^2 = \|x\|^2$ and $T_\theta$ is a contraction,
we see that
\[
\|x\| \geq \|T_\theta \Gamma x\| \geq
\|\Gamma\Gamma^* T_\theta \Gamma x\|=\|\Gamma\varphi(A^*) x\| = \|\varphi(A^*) x \| =\|x\|.
\]
Therefore we have equality, and thus,
\[
T_\theta \Gamma x = \Gamma\Gamma^* T_\theta \Gamma x = \Gamma\varphi(A^*) x
\qquad (\mbox{when } 0\neq x \in \mathfrak{D}_{\varphi(A^*)}^\perp).
\]
By taking the Fourier transform of both sides, and using \eqref{clt}, we obtain
\[
\theta(z) = \frac{C(I-z A)^{-1} \varphi(A^*) x}{C(I-z A)^{-1}x} = \varphi(z).
\]
So $\theta(z) = \varphi(z)$. In other words, if
$ \theta $ is a function in $H^\infty$ such that
$\widetilde{\theta}(A) = \widetilde{\varphi}(A)$ and $\|\theta\|_\infty \leq 1$,
then $\theta(z) = \varphi(z)$. Replacing $\widetilde{\varphi}$ by $\varphi$
shows that if $ \theta $ is a function in $H^\infty$ such that
$ \theta (A) =  \varphi (A)$ and $\|\theta\|_\infty \leq 1$,
then $\widetilde{\theta}(z) = \widetilde{\varphi}(z)$ and therefore $\theta(z) = \varphi(z)$.
\epr

\bigskip

Finally, it is noted that this result is an application    of
the Sz.-Nagy-Foias commutant lifting theorem; see Corollary 2.7
page 142 of \cite{FFGK} and is also deeply connected to
some of the results in  \cite{AAK}.

 \begin{proposition}\label{prop-phi}
 Let $A$ be a stable contraction on a finite dimensional space
 $\mathcal{X}$ and $\varphi$ a Blaschke product. If
 $\deg(\varphi) \geq \dim(\mathcal{X})$, then
 the defect index for $\varphi(A)$ equals $\dim(\mathcal{X})$,
 or equivalently,
 $\varphi(A)^*\varphi(A)$ has no eigenvalues on the unit circle.
 \end{proposition}

\noindent {\sc Proof.}
 Let $C$ be any operator from $\mathcal{X}$ onto $\mathcal{E}$ such
 that $I -A^*A = C^*C$. Then the observability operator
 $\Gamma $ formed by $\{C,A\}$ is an isometry
 from $\mathcal{X}$ into $\ell_+^2(\mathcal{E})$. Moreover, the range of
 $\Gamma $ is an invariant subspace for the backward
 shift $S^*$ on $\ell_+^2(\mathcal{E})$.
 Recall that if $u(z)$ is any function in $H^\infty$, then
 $\|u(A)\| \leq \|u\|_\infty$; see the functional calculus in \cite{sznf}.
 Because  $\varphi$ is an inner function, $\widetilde{\varphi}(A)$ is a contraction.
  For $x$ in $\mathcal{X}$, we have
 \begin{align*}
 \|\big(I - \widetilde{\varphi}(A)^*\widetilde{\varphi}(A)\big)^{\frac{1}{2}}x\|^2
 &= \langle\big(I - \widetilde{\varphi}(A)^*\widetilde{\varphi}(A)\big)x,x\rangle  =
 \|x\|^2 - \|\widetilde{\varphi}(A) x\|^2 \\
 &= \|\Gamma x\|^2 - \|\widetilde{\varphi}(A) x\|^2 =
 \|\Gamma x\|^2 - \|T_{\varphi I} T_{\varphi I}^* \Gamma x\|^2\\
 &= \|(I -  T_{\varphi I}T_{\varphi I}^*) \Gamma x\|^2=
 \|P_{_{\mathfrak{H}(\varphi I)}}\Gamma x\|^2.
 \end{align*}
 Hence there exists a unitary operator
 $\Psi$ from the
 range of $P_{_{\mathfrak{H}(\varphi I)}}\Gamma $
 onto the range of the defect operator
 $\big(I - \widetilde{\varphi}(A)^*\widetilde{\varphi}(A)\big)^{\frac{1}{2}}$
 such that
 \begin{equation}\label{defect}
 \Psi P_{_{\mathfrak{H}(\varphi I)}}\Gamma  =
 \big(I - \widetilde{\varphi}(A)^*\widetilde{\varphi}(A)\big)^{\frac{1}{2}}.
 \end{equation}
 In particular, the defect index $\mathfrak{d}_{_{\widetilde{\varphi}(A)}}$
 equals the rank of $P_{_{\mathfrak{H}(\varphi I)}}\Gamma $.

 Notice that $x$ is in the kernel of $P_{_{\mathfrak{H}(\varphi I)}}\Gamma $
 if and only if $x$ is a vector with  eigenvalue $1$ for
 $\widetilde{\varphi}(A)^*\widetilde{\varphi}(A)$. In this case,
 \[
 0 = P_{_{\mathfrak{H}(\varphi I)}}\Gamma  x =
 \Gamma  x - T_{\varphi I}T_{\varphi I}^* \Gamma  x.
 \]
 Hence
 \[
 \Gamma x =  T_{\varphi I} \Gamma  \widetilde{\varphi}(A)x
 \qquad (\mbox{for } x \in \mathfrak{D}_{\widetilde{\varphi}(A)}^\perp).
 \]
 By taking the Fourier transform
 \begin{equation}  
 C(I-z A)^{-1}x = \varphi(z) C(I-z A)^{-1}\widetilde{\varphi}(A)x
 \qquad (\mbox{for } x \in \mathfrak{D}_{\widetilde{\varphi}(A)}^\perp).
 \end{equation}
 This means that the numerator of $C(I-z A)^{-1}x$ has at least
 $\deg(\varphi)$ zeros. However, the numerator of $C(I-z A)^{-1}x$
 is a polynomial of degree at most $\dim(\mathcal{X}) -1 < \deg(\varphi)$.
 Therefore $x=0$, and
 $\widetilde{\varphi}(A)\widetilde{\varphi}(A)^*$
has no eigenvalue on the unit circle. The same argument
applies by replacing $\widetilde{\varphi}(z)$ by
$\varphi(z)$.  This completes the proof.
\epr

 \paragraph{Winding numbers in the scalar case}
Let $R$ be the rigid function in $L^\infty$ defined by
$R(z) = \varphi(z) \overline{m(1/\bar{z})} = (\varphi m^*)(z)$,
where both $\varphi$ and $m$ are Blaschke products.
The {\it winding number} for $R$ is given by $\deg(\varphi) - \deg(m)$.
As before, let $Y$ be the contraction from
$\mathfrak{H}(m)$ into $\mathfrak{H}(\varphi)$ given by
$Y = P_{_{\mathfrak{H}(\varphi)}}\vert \mathfrak{H}(m)$.
By consulting Lemma \ref{lem-mphi}, we see that
\begin{enumerate}
  \item The winding number for $R = \varphi m^*$ equals
  $\dim\big(\im(Y)^\perp\big)$ when\\ $\deg(m) \leq \deg(\varphi)$.
  \item The winding number for $R = \varphi m^*$ equals
  $-\dim\big(\kr(Y)\big)$ when\\ $\deg(\varphi) \leq \deg(m)$.
  \item Hence  the winding number for  $R = \varphi m^*$ equals $-\ind(Y)$
 where $\ind(Y)$ denotes  the Fredholm index
 $\dim(\kr(Y)) - \dim(\kr(Y^*))$)  for $Y$.
\end{enumerate}

Recall that the Toeplitz operator $T_R$ admits a  decomposition of the form:
\[
T_R = T_\varphi T_m^* + H_\varphi Y H_m^*;
\]
see \eqref{svdTR}.
Here $Y$ is the finite dimensional contraction from $\mathfrak{H}(\widetilde{m})$ into
$\mathfrak{H}(\widetilde{\varphi})$ determined  by
$Y = P_{_{\mathfrak{H}(\widetilde{\varphi})}}\vert \mathfrak{H}(\widetilde{\varphi})$.
Hence
$T_R$ is Fredholm. Moreover, $T_R$ and $Y$ have the
same Fredholm index.   Hence
\begin{align}\label{frecind}
\ind(T_R) &= \dim(\kr (T_R)) - \dim(\kr(T_R^*))
 = \ind(Y) =\deg(m) - \deg(\varphi).
\end{align}
Finally, the Fredholm index of $T_R$ equals minus the
winding number of $R = \varphi m^*$.


\setcounter{equation}{0}
\section{Point evaluation and multidimensional systems}

Let us return to the multidimensional case
and provide a proof of one of our main results, Theorem \ref{mainthm00}.

As before, assume that $V$ and $W$ are two bi-inner rational functions
in $H^\infty(\mathcal{E},\mathcal{E})$.
Let $\{A_v \mbox{ on } \mathcal{X}_v, B_v,C_v,D_v\}$
and $\{A_w \mbox{ on } \mathcal{X}_w, B_w,C_w,D_w\}$
be stable unitary realizations of $V$ and $W$, respectively.
Recall that  $R=VW^*$ and the Toeplitz operator  $T_R$ on $\ell_+^2(\mathcal{E})$ admits
a decomposition of the form:
\begin{equation}\label{TRdef8}
T_R = T_V T_{_W}^* + H_V Y H_{W}^*.
\end{equation}
Here $Y$ is the contraction defined by
\begin{equation}\label{Ydef8}
Y = P_{_{\mathfrak{H}(\widetilde{V})}}\vert \mathfrak{H}(\widetilde{W})
: \mathfrak{H}(\widetilde{W}) \rightarrow \mathfrak{H}(\widetilde{V}).
\end{equation}
Recall that Proposition \ref{prop-R} shows that the dimensions of $\kr(Y)$ and
$ \kr (T_R) $ are equal.
Therefore we are interested in calculating the dimension of $\kr(Y)$ in terms
of the realizations of $V$ and $W$.
Since $H_W = \Gamma_{w} \Upsilon_w$, and the
range of the  Hankel operator $\im(H_W^*) = \mathfrak{H}(\widetilde{W})$, we have
\[
\mathfrak{H}(\widetilde{W} ) =  \im ( \Upsilon_w^*)
\quad \mbox{where}\quad
 \Upsilon_w^* = \begin{bmatrix}
 B_w^* \\  B_w^* A_w^* \\  B_w^* A_w^{*2} \\  \vdots \\
 \end{bmatrix}: \mathcal{X}_w \rightarrow \ell_+^2(\mathcal{E}).
\]
Furthermore, $\Upsilon_w $ is a co-isometry and $ A_w^* $
 is a stable contraction on $\mathcal{X}_w$
satisfying $I = A_w A_w^* + B_w B_w^* $, where $B_w$ is an operator mapping
$\mathcal{E}$ into $\mathcal{X}_w$.
Recall that if $ V(z) = \sum_{0}^\infty z^k V_k$ is the Taylor series
expansion for $V(z)$, then
\[
V_0 = D_v \qquad\mbox{and}\qquad
V_k = C_v A_v^{k-1}B_v \quad (\mbox{for } k \geq 1).
\]
Let $C_\circ$ be the operator mapping $\mathcal{X}_w$ into $\mathcal{E}$ defined by
\begin{equation}\label{defCo}
C_\circ  =   \sum_{k=0}^\infty V_k B_w^* (A_w^*)^k
= D_v B_w^* + C_v\left(\sum_{k=0}^\infty  A_v^{k} B_v B_w^* A_w^{*k}\right) A_w^*.
\end{equation}
Because the operators $A_v$ and $A_w$ are both stable,
the operator $C_\circ$ is well defined.
In fact, $C_\circ$ can be computed by first solving the Stein
equation
\begin{equation}\label{lyapom}
\Omega = A_v \Omega A_w^* + B_v B_w^*.
\end{equation}
Because $A_v$ and $A_w$ are stable, the solution $\Omega$ to this
Stein equation is unique and given by
\begin{equation} \label{defOmega}
\Omega = \sum_{k=0}^\infty  A_v^{k} B_v B_w^* A_w^{*k}.
\end{equation}
Therefore
\begin{equation}\label{Ccirc}
C_\circ = D_v B_w^* + C_v \Omega A_w^*.
\end{equation}

Now let $\Gamma_\circ$ be the operator
from $\mathcal{X}_w$ into $\ell_+^2(\mathcal{E})$
defined by  $ \Gamma_\circ = T_{\widetilde{V}}^* \Upsilon_w^* $.
Then a simple calculation shows that
\begin{equation}\label{Gammacirc}
 \Gamma_\circ  =    T_{\widetilde{V}}^* \Upsilon_w^*  =
\begin{bmatrix} C_\circ \\
              C_\circ A_w^* \\
              C_\circ A_w^{*2} \\
              \vdots \\
\end{bmatrix}: \mathcal{X}_w  \rightarrow \ell_+^2(\mathcal{E}).
\end{equation}
Notice  that $\Gamma_\circ$ is the observability operator
determined by the pair $\{C_\circ,A_w^*\}$.

It is emphasized that if $V = \varphi I$ where $\varphi$
is an inner function in $H^\infty$,
then $C_\circ = B_w^* \varphi(A_w^*)$. In other words, in the scalar case
\begin{equation}\label{qqq}
\Gamma_\circ =
T_{\widetilde{\varphi}I}^* \Upsilon_w^* = \Upsilon_w^* \varphi(A_w^*)
\end{equation}
and the operator $C_\circ$ plays the role of   $B_w^*\varphi(A_w^*)$
in the multidimensional case. Note that in the general case,
$C_\circ$ is the left point evaluation
of $ A_w^* $ with respect to $  B_w^* $.
For a further discussion on multidimensional
function evaluation with applications to $H^\infty$ interpolation theory,
see Section 1.2 page 15 of \cite{FFGK}.

The next lemma together with $\mathfrak{n}(T_R) = \mathfrak{n}(Y)$ provides the proof of Part 1 and the first half of Part 2
of Theorem \ref{mainthm00}.

\begin{lemma}
Let $V$ and $W$ be  bi-inner rational  functions  in
$H^\infty(\mathcal{E},\mathcal{E})$.
Consider the contraction mapping $Y$ on $\mathfrak{H}(\widetilde{W})$
into $\mathfrak{H}(\widetilde{V})$  determined  by
\begin{equation}\label{Y43}
 Y = P_{_{\mathfrak{H}(\widetilde{V})}}\vert\mathfrak{H}(\widetilde{W}):
 \mathfrak{H}(\widetilde{W}) \rightarrow \mathfrak{H}(\widetilde{V}).
\end{equation}
Let $ C_\circ$ be the operator from $\mathcal{X}_w$ into $\mathcal{E}$
defined by \eqref{defCo} or \eqref{lyapom} and \eqref{Ccirc}.
Let $ Q$ be the unique solution to the Stein equation
\begin{equation}\label{defQ}
 Q = A_w  Q A_w^* + C_\circ^* C_\circ.
\end{equation}
Then $Q$ is a positive contraction. Moreover, $Y$ and $I-Q$ have the same rank.
In particular,
\begin{equation}\label{eig1}
\kr(Y) = \Upsilon_w^*\kr(I - Q)
\quad \mbox{and}\quad \dim(\kr(Y))= \dim(\kr( I- Q)).
\end{equation}
Furthermore, we have
\begin{equation}\label{eig1dual}
\dim(\kr(Y^*))   = \dim(\mathfrak{H}(\widetilde{V})) - {\rm rank} (I -  Q).
\end{equation}
 If   $V = \varphi I$ where $\varphi$ is a Blaschke product,
 then $Q = \varphi(A_w^*)^*\varphi(A_w^*)$.
\end{lemma}

\bpr
Notice  that $ Q =  \Gamma_\circ^*  \Gamma_\circ $ is
the unique solution to the Stein equation
\[
 Q = A_w  Q A_w^* +  C_\circ^* C_\circ.
\]
Recall that the orthogonal projection
  $ P_{_{\mathfrak{H}(\widetilde{V})}}
= I - T_{\widetilde{V}} T_{\widetilde{V}}^* $.
Moreover,  $\mathfrak{H}(\widetilde{W})$ equals the
range of the isometry $\Upsilon_w^*$.
So $h$ is in $\mathfrak{H}(\widetilde{W})$
if and only if $h =  \Upsilon_w^* x$ for some
$x$ in  $\mathcal{X}_w$. In fact, $x = \Upsilon_w h$. Now observe that
\begin{align*}
  \|P_{_{\mathfrak{H}(\widetilde{V})}} \Upsilon_w^* x \|^2 &=
  \|\big(I - T_{\widetilde{V}} T_{\widetilde{V}}^* \big) \Upsilon_w^* x\|^2
  =
  \| \Upsilon_w^* x\|^2 -
  \| T_{\widetilde{V}} T_{\widetilde{V}}^* \Upsilon_w^* x \|^2\\
&=
\| \Upsilon_w^* x\|^2 - \|  T_{\widetilde{V}}^*\Upsilon_w^* x\|^2
 = \|x\|^2 - \|  \Gamma_\circ  x\|^2 = \|x\|^2 - \langle Q x,x\rangle.
\end{align*}
Hence $Q$ is a positive contraction, and thus,
\[
 \|P_{_{\mathfrak{H}(\widetilde{V})}} \Upsilon_w^* x \|^2
  = \|(I -  Q)^\frac{1}{2}x\|^2 \qquad (x\in \mathcal{X}).
\]
So there exists a unitary operator $\Psi$ mapping
the range of $P_{\mathfrak{H}(\widetilde{V} )} \Upsilon_w^* $ onto
the range $(I -  Q)^\frac{1}{2}$ such that
\begin{equation}\label{defect6}
\Psi Y  \Upsilon_w^* x =
\Psi P_{\mathfrak{H}(\widetilde{V} )}  \Upsilon_w^* x
= \big(I -  Q\big)^\frac{1}{2} x
\qquad (\mbox{for } x\in  \mathcal{X}_w ).
\end{equation}
In particular, $Y$ and $I-Q$ have the same rank.
Recall that $\mathfrak{H}( \widetilde{W} )$ equals the range of
$\Upsilon_w^*$
and that $Y = P_{_{\mathfrak{H}(\widetilde{V} )}} \vert  \mathfrak{H}(\widetilde{W}) $.
Therefore
\[
\kr(Y) = \Upsilon_w^*\kr(I - Q)
\quad \mbox{and}\quad
\dim(\kr(Y))= \dim (\kr( I- Q)) .
\]
In other words, $\dim (\kr(Y) )$ equals the number of eigenvalues of $ Q$ equal to $1$ counting multiplicities.

By applying the previous result to $Y^*$, we obtain
\[
\dim (\kr(Y^*)) = \dim (\mathfrak{H}( \widetilde{V} ))  - \rank(I -  Q).
\]

Finally, it is noted that if $V = \varphi I$ where
$\varphi$ is an inner function in $H^\infty$, then
\[
\Gamma_\circ = \Upsilon_w^*\varphi(A_w^*);
\]
see \eqref{qqq}. Because $\Upsilon_w $ is a co-isometry, we have
\[
Q = \Gamma_\circ^*\Gamma_\circ =
\varphi(A_w^*)^* \Upsilon_w \Upsilon_w^* \varphi(A_w^*)
= \varphi(A_w^*)^* \varphi(A_w^*).
\]
\epr

The following lemma implies Part 2 of Theorem \ref{mainthm00}.

\begin{lemma}
Let $V$ and $W$ be two rational bi-inner functions in
$H^\infty(\mathcal{E},\mathcal{E})$, and  $\psi$ be an
inner function in $H^\infty$.
Let $Y_\psi$ be the
contraction from $\mathfrak{H}(\widetilde{W})$ into
$\mathfrak{H}(\psi \widetilde{V})$ defined by
$Y_\psi = P_{\mathfrak{H}( \psi  \widetilde{V} )}
\vert {\mathfrak{H}( \widetilde{W} ) }$.
Then $\dim (\kr(Y_\psi )) $ is equal to the multiplicity of $1$ as an eigenvalue of
$\psi(A_w) Q \psi(A_w)^* $.
Also
\[
\dim (\kr(Y_\psi^*)) = \deg(\psi) \dim \mathcal{E} + \dim (\mathfrak{H}( \widetilde{V} )) -
\rank (I - \psi(A_w) Q \psi(A_w)^* ).
\]
In particular, if $\psi(z) = z^k$, then
 $\dim (\kr(Y_{z^k})) $ is equal to the multiplicity of $1$ as an eigenvalue of
$ A_w^k Q A_w^{*k} $.
Also
\[
\dim (\kr(Y_{z^k}^*))  = k\dim (\mathcal{E} )+ \dim (\mathfrak{H}( \widetilde{V} )) -
\rank (I -  A_w^k Q  A_w^{*k} ).
\]
\end{lemma}

\bpr
Set $\Gamma = \Upsilon_w^*$.
Then for $x$ in $\mathcal{X}_w$ we have
\begin{equation*}
\Vert P_{\mathfrak{H}( \psi \widetilde{V}) }  \Gamma x \|^2 =
\| ( I - T_{  \psi\widetilde{ V}} T_{ \psi \widetilde{  V}}^* ) \Gamma x  \|^2 =
\| \Gamma x \|^2 - \| T_{\psi \widetilde{V}} T_{\psi \widetilde{V}}^* \Gamma x\|^2 .
\end{equation*}
Since  $ T_{ \psi\widetilde{ V}} $ is an isometry, and using \eqref{Gammacirc}, we have
\begin{align*}
\Vert P_{\mathfrak{H}(  \psi\widetilde{ V}) }  \Gamma x \|^2 &=
\| \Gamma x \|^2 - \| T_{ \psi \widetilde{  V}}^* \Gamma x\|^2
 = \| \Gamma x \|^2 - \| T_{  \widetilde{ V}}^* T_{\psi }^*\Gamma x\|^2\\
&= \| \Gamma x \|^2 - \| T_{  \widetilde{ V}}^*\Gamma \psi(A_w)^*x\|^2\\
&= \| \Gamma x \|^2 - \| \Gamma_\circ \psi(A_w)^*x\|^2.
\end{align*}
Recall that $ \Gamma_\circ^*  \Gamma_\circ  =  Q $.
Therefore
\begin{align*}
\Vert P_{\mathfrak{H}( \psi \widetilde{V}) }  \Gamma x \|^2 &=
\|x\|^2 - \langle \psi(A_w) Q \psi(A_w)^* x,x\rangle\\
&=
\|(I - \psi(A_w) Q \psi(A_w)^*)^\frac{1}{2}x\|^2.
\end{align*}
So there exists a unitary operator $\Psi$ mapping
the range of $P_{\mathfrak{H}( \psi \widetilde{V}) }  \Gamma $ onto
the range of
$(I - \psi(A_w) Q \psi(A_w)^*)^\frac{1}{2} $ such that
\[
\Psi P_{\mathfrak{H}( \psi \widetilde{V}) }  \Gamma x =
\big(I - \psi(A_w) Q \psi(A_w)^*\big)^\frac{1}{2} x
\qquad (\mbox{for } x \in  \mathcal{X}_w ).
\]
Recall that $\mathfrak{H}( \widetilde{W} ) = \Gamma \mathcal{X}_w$.
Therefore, when
$Y_\psi = P_{{\mathfrak{H}(\psi \widetilde{V} )}} \vert {\mathfrak{H}( \widetilde{W} )}$, we have
\[
\dim \left( \kr \left( P_{\mathfrak{H}(\psi \widetilde{V} )}
\vert {\mathfrak{H}( \widetilde{W})} \right)\right)
= \dim\left( \kr \left( I- \psi(A_w) Q \psi(A_w)^* \right) \right).
\]
In other words, the dimension of $ \kr(Y_\psi )$ is equal to the multiplicity of $1$ as an eigenvalue
of $ \psi(A_w) Q \psi(A_w)^* $.

Using $\mathfrak{H}(\psi \widetilde{V} )  = \mathfrak{H}(\psi I_\mathcal{E})\oplus
T_\psi \mathfrak{H}( \widetilde{V})$,
we also have
\[
\dim \left(\kr(Y_\psi^* ) \right) = \deg(\psi) \dim( \mathcal{E}) +
\dim (\mathfrak{H}( \widetilde{V} )) -
\mbox{rank}(I - \psi(A_w) Q \psi(A_w)^* ).
\]
\epr

The previous lemma combined with $\mathfrak{n}(T_{z^kR}) = \mathfrak{n}(Y_{z^k})$, yields
\[
\mathfrak{n}(T_{z^k R})  = \dim \left( \kr (I - A_w^{k} Q A_w^{*k} )\right).
\]
This completes the proof of Part 2 of Theorem \ref{mainthm00}.

Since Part 3 of Theorem \ref{mainthm00} follows from the Parts 1 and 2 the proof of Theorem \ref{mainthm00} is complete now.

The next Corollary provides similar formulas for the positive Wiener-Hopf indices.
First we define $C_{\circ\ast}$ as follows.
Solve the Stein
equation
\begin{equation}\label{lyapom_ast}
\Omega_\ast = A_w \Omega_\ast A_v^* + B_w B_v^*.
\end{equation}
(Note that $\Omega_*=\Omega^*$.)
Because $A_v$ and $A_w$ are stable the solution $\Omega_*$ to this
Stein equation is unique and given by
\begin{equation} 
\Omega_\ast = \sum_{k=0}^\infty  A_w^{k} B_w B_v^* A_v^{*k}.
\end{equation}
Put
\begin{equation}\label{Ccircd}
C_{\circ\ast} = D_w B_v^* + C_w \Omega_\ast A_v^*.
\end{equation}
Now let $ Q_\ast $ be the unique solution of the Stein equation
\begin{equation}\label{defQast}
Q_\ast = A_v Q_\ast A_v^\ast + C_{\circ\ast}^* C_{\circ\ast} .
\end{equation}

\begin{corollary}\label{maincol}
Let $V$ and $ W $ be given by \eqref{defV00} and \eqref{defW00} and $ R = V W^* $.
Furthermore let $ Q_\ast $ be defined by \eqref{defQast}.
Put
\begin{equation}\label{defnuk}
\nu_k = \dim (\kr(I - A_v^{k-1} Q_\ast A_v^{*(k-1)} ))
- \dim ( \kr (I - A_v^{k} Q_\ast A_v^{*k} ) ).
\end{equation}
Then the positive Wiener-Hopf indices $ \omega_1 , \ldots , \omega_q $
of $ T_R $ are given by
\[
\omega_j = \# \{ k \, : \, \nu_k \geq j \}, \quad (j = 1 , \ldots, q = \nu_1 ).
\]
\end{corollary}

\bpr
Notice $ R^* = W V^* $.
Recall $ R^*(z) $ is defined by $R^*(z) = (R(\frac{1}{\overline{z}} )^* $.
If
\[
R(z) = W_-(z) \ \diag( z^{\kappa_j} )_{j=1}^m \ W_+(z) ,
\]
with  $ W_+ $ and its inverse are analytic on $ \overline{\BD } $ and
$W_- $ and its inverse are analytic on the complement of $ \BD $,
 then
\[
 R^*(z) = W_+^*(z) \diag( z^{-\kappa_j} )_{j=1}^m W_-^*(z) .
\]
Moreover $ W_-^* $ and its inverse are analytic on $  \overline{\BD } $ and
$W_+^* $ and its inverse are analytic on the complement of $ \BD $
in $\BC \cup \infty $.
This shows that $ - \kappa_1, \ldots -\kappa_m $
are the Wiener-Hopf indices of $ R^* $.
The positive Wiener-Hopf indices of $ R $ are the opposite to the
negative Wiener-Hopf indices of $ R^* $.
So the Corollary is immediate from applying Theorem \ref{mainthm00} to $ R^* $.
\epr

\bigskip\noindent

\subsection{An example}
As an illustration of Theorem \ref{mainthm00} we present the following example
that also appears on page 706 in \cite{GKRa}. To this end, let
\[
R(z) = \begin{bmatrix}
z^{-4}&0&0&0&0\\  0&z^{-2}&0&0&0\\  0&0&1&0&0\\  0&0&0&z^3&0\\  0&0&0&0&z^5
\end{bmatrix}.
\]
Then $ R(z) $ factors as $ R(z) = V(z) W^*(z)$, where
\[
V(z) = \begin{bmatrix}
1&0&0&0&0\\  0&1&0&0&0\\  0&0&1&0&0\\  0&0&0&z^3&0\\  0&0&0&0&z^5
\end{bmatrix}
\qquad
W(z) = \begin{bmatrix}
z^{4}&0&0&0&0\\  0&z^{2}&0&0&0\\  0&0&1&0&0\\  0&0&0&1&0\\  0&0&0&0&1
\end{bmatrix} .
\]
We shall use the following notations: $ J_n(0) $ denotes the standard upper
triangular Jordan block with eigenvalue zero of size $ n \times n$,
and $e_j $ denotes
the $j$'th standard unit vector in a Euclidean space, with a one in the $j$-th
position and zeros everywhere else. Note that the size of $e_j $ depends on the
particular choice of the Euclidean space. Finally, it is noted that a
stable unitary realization of $z^n$ is given by
\[
z^n = z e_1(I-z J_n(0))^{-1}e_n.
\]

Motivated by the previous realization, the factors $ V(z) $ and $ W(z) $
can be given by the following stable unitary realizations:
\[
V(z) = D_v + z C_v ( I_{\BC^8} - z A_v )^{-1} B_v ,
\]
where
\begin{align*}
& A_v = J_3(0) \oplus J_5(0), \qquad\qquad\quad\,
B_v = \begin{bmatrix} 0 & 0 & 0 & e_3 & e_8 \end{bmatrix} , \\
& C_v = \begin{bmatrix} 0 & 0 & 0 & e_1 & e_4 \end{bmatrix}^T , \qquad
D_v = \begin{bmatrix} e_1 & e_2 & e_3 & 0 & 0 \end{bmatrix} ,
\end{align*}
and
\[
W(z) = D_w + z C_w ( I_{\BC^6} - z A_w )^{-1} B_w ,
\]
where
\begin{align*}
& A_w = J_4(0) \oplus J_2(0), \qquad \qquad\quad\,
B_w = \begin{bmatrix} e_4 & e_6 & 0 & 0 & 0 \end{bmatrix} ,\\
& C_w = \begin{bmatrix} e_1& e_5 & 0  & 0 & 0 \end{bmatrix}^T , \qquad
D_w = \begin{bmatrix} 0 & 0 & e_3 & e_4 & e_5  \end{bmatrix} .
\end{align*}

We choose to calculate $ Q $ according to the definitions.
We determine $ C_\circ $ from
\[
C_\circ = \sum_{k=0}^\infty V_k B_w^* ( A_w^*)^k ,
\]
where $V_k $ is the $k$-th  Fourier coefficient of $ V $.
Then
\[
V_0 = D_v , \quad V_3 = \begin{bmatrix} 0 & 0 & 0 & e_4 & 0  \end{bmatrix}, \quad
V_5 = \begin{bmatrix} 0 & 0 & 0 & 0 & e_5 \end{bmatrix}
\]
and for all other $ k $ one has  $ V_k = 0$.
Also notice that $ ( A_w^*)^k = 0 $ whenever $ k\geq 4 $.
Therefore we have, also using $ V_3 B_w^* = 0 $, that
\[
C_\circ = V_0 B_w^* + V_3 B_w^* ( A_w^*)^3 = V_0 B_w^* = B_w^*.
\]
Next we compute $ \Gamma_\circ $.
Since $ C_\circ= B_w^* $, we have that $ \Gamma_\circ = \Upsilon_w^* $.
Now use that the realization of $W$ is unitary to conclude that
\[
Q = \Gamma_\circ^* \Gamma_\circ = \Upsilon_w \Upsilon_w^* = I .
\]

So, we get that
\begin{align*}
& \dim (\kr (I-Q)) = 6, \quad \dim(\kr(I- A_w  Q A_w^* ) )= 4, \\
& \dim\left(\kr \left(I - A_w^2 Q (A_w^*)^2 \right) \right)= 2 , \quad
\dim\left(\kr \left(I - A_w^3 Q (A_w^*)^3 \right)\right) = 1, \\
& \dim\left(\kr \left(I - A_w^k Q (A_w^*)^k \right) \right)= 0 \mbox{\ \ for  } k \geq 4.
\end{align*}
Finally, we see that
\[
\mu_1 = 6-4=2, \quad \mu_2 = 4-2=2,\quad \mu_3 = 2-1=1 \mbox{ and } \mu_4=1-0=1.
\]
Using this we have
\[
\kappa_1 = \#\{k:\mu_k \geq 1\} = 4
\quad \mbox{and}\quad \kappa_2 = \#\{ k:\mu_k \geq 2 \} = 2.
\]
Therefore the negative Wiener-Hopf indices for $R$ are $\{-4,-2\}$.

Because this example is simple, we could compute $Q$ by hand.
However, in general one has to compute the solutions to the
corresponding  Stein equations to determine $Q$.


\setcounter{equation}{0}
\section{Direct connection with earlier results}
In this section,  we will  establish a direct  connection with the formulas
in Theorem 3.3 in the paper \cite{GKRb}. Let us begin by presenting  this Theorem.
To accomplish  this we need some notation.
As before, let $R = VW^*$ where $V$ and $W$ are two rational
bi-inner functions in $H^\infty(\mathcal{E},\mathcal{E})$
with corresponding stable unitary realizations
 $\{A_v \mbox{ on } \mathcal{X}_v, B_v,C_v,D_v \}$
and $\{A_w \mbox{ on } \mathcal{X}_w, B_w,C_w,  D_w\}$.

Recall that $Y$ is the operator mapping $\mathfrak{H}(\widetilde{W})$
into $\mathfrak{H}(\widetilde{V})$ determined by
$Y = P_{_{\mathfrak{H}(\widetilde{V})}}\vert \mathfrak{H}(\widetilde{W})$.
Moreover,  $\Upsilon_v^*$ is an isometry from  $\mathcal{X}_v$
into $\ell_+^2(\mathcal{E})$ whose range equals $\mathfrak{H}(\widetilde{V})$, and
$\Upsilon_w^*$ is an isometry from  $\mathcal{X}_w$
into $\ell_+^2(\mathcal{E})$ whose range equals $\mathfrak{H}(\widetilde{W})$.
An important role in \cite{GKRb} is played by the operator
\begin{equation}\label{defX}
X = \Upsilon_v \Upsilon_w^*:
\mathcal{X}_w \rightarrow  \mathcal{X}_v.
\end{equation}
It is emphasized that $X$ is the unique solution of the Stein equation:
\begin{equation}\label{Stein}
X = A_v X A_w^* + B_v B_w^*.
\end{equation}
(Observe that by \eqref{lyapom} we have $X=\Omega$, but we will use $X$ here
to explain the connection with \cite{GKRb}.)
Using the fact that the range of $\Upsilon_w^*$ equals
$\mathfrak{H}(\widetilde{W})$
and the  $\kr\big(\Upsilon_v \big)^\perp$ equals $\mathfrak{H}(\widetilde{V})$,
we see that
\begin{equation}\label{defX1}
X = \Upsilon_v Y  \Upsilon_w^*: \mathcal{X}_w \rightarrow  \mathcal{X}_v.
\end{equation}
Consider the pair $\{C,A\}$ where $A$ is an operator on $\mathcal{X}$ and
$C$ maps $\mathcal{X}$ into $\mathcal{E}$. Then let us set
\[
\kr_{m}  ( C,A) = \bigcap_{j=0}^{m} \kr \left( CA^j \right).
\]
We now are ready to rephrase Theorem 3.3 in the paper \cite{GKRb}.

\begin{theorem}\label{Thm33}
The number $ s$ of negative Wiener Hopf indices
of the function $ R = V W^* $ is given by
\[
s =  \dim (\kr X )- \dim \left(\kr \begin{bmatrix}  B_w^* \\ X A_w^* \end{bmatrix}\right).
\]
Let
\[
\nu_k = \dim\left( \kr_{k-1}
\left( \begin{bmatrix}  B_w^* \\ X A_w^* \end{bmatrix} , A_w^*\right)\right) -
\dim\left( \kr_{k}
\left( \begin{bmatrix}  B_w^* \\ X A_w^* \end{bmatrix} , A_w^*\right)\right).
\]
Then the negative Wiener-Hopf indices $ -\kappa_1 , \ldots , -\kappa_p $
of $ T_R $ are given by
\[
\kappa_j = \# \{ k \, : \, \nu_k \geq j \}, \quad (j = 1 , \ldots, p = \nu_1 ).
\]
\end{theorem}

\bigskip
Comparing Theorem \ref{mainthm00} with Theorem \ref{Thm33} reveals that apparently
\[
\dim\left( \kr_{k}
\left( \begin{bmatrix}  B_w^* \\ X A_w^* \end{bmatrix} , A_w^*\right)\right)
= \dim \left( \kr (I - A_w^{k} Q A_w^{*k} ) \right)
\]
We will show this directly, that is not using the detour via Wiener Hopf indices.
It will prove that Theorem \ref{Thm33} is equivalent to Theorem \ref{mainthm00}.

The first step is the following lemma.

\begin{lemma}
The operator $X$ satisfies the identities:
\begin{equation}\label{threeX}
X = \Gamma_v^* T_R \Gamma_w=
 \Upsilon_v Y  \Upsilon_w^*
 \quad \mbox{and} \quad Q = I-X^*X.
\end{equation}
\end{lemma}

\bpr
Since $ \ell^2_+( \mathcal{E}) = \mathfrak{H}(V) \oplus \im(T_V)$ and
(see equality \eqref{GammaH})
$P_{\mathfrak{H}(V) } = \Gamma_v \Gamma_v^* $ we have that
$ \Gamma_v \Gamma_v^* T_V = 0 $.
Now use that $ \Gamma_v $ is an isometry to see that $ \Gamma_v^* T_V = 0 $.
Recall that $R = V W^*$ and $T_R = T_V T_W^* + H_V H_W^* $.
Using this with $H_V = \Gamma_v \Upsilon_v$ and $H_W = \Gamma_w \Upsilon_w$, we have
\begin{equation}\label{DefX6}
 \Gamma_v^* T_R \Gamma_{w}=
\Gamma_v^* \Gamma_v \Upsilon_v \Upsilon_w^* \Gamma_w^* \Gamma_w =
 \Upsilon_v \Upsilon_w^* = X.
\end{equation}
This yields the first equality in \eqref{threeX}.

According to \eqref{defect6}, we have
$\Psi Y  \Upsilon_w^* = \big(I -  Q\big)^\frac{1}{2}$
where $\Psi$ is a unitary operator
from the range of $Y  \Upsilon_w^*$ onto the range of
$\big(I -  Q\big)^\frac{1}{2}$. Now observe that
\begin{align*}
X^*X & = \Upsilon_w Y^* \Upsilon_v^* \Upsilon_v Y  \Upsilon_w^*
= \Upsilon_w Y^* P_{_{\mathfrak{H}(\widetilde{V})}} Y  \Upsilon_w^* \\
&= \Upsilon_w  Y^*  Y  \Upsilon_w^*=
\Upsilon_w Y^* \Psi^* \Psi Y  \Upsilon_w^*\\
&=I-Q.
\end{align*}
Therefore $I-X^*X = Q$.
\epr

For further results on the  operator $X$ see  \cite{GKRb} and \cite{FK}.

\begin{lemma}
For $ k = 1,2,3,\cdots $ we have
\begin{equation}\label{kerk}
\kr \left( I - \left( A_w^\ast \right)^k Q A_w^k \right) = \kr
\left[\begin{array}{l}
  \begin{bmatrix} B_w^\ast \\ X A_w^\ast \end{bmatrix}\\\\
  \begin{bmatrix} B_w^\ast \\ X A_w^\ast \end{bmatrix} A_w^*\\\\
  \begin{bmatrix} B_w^\ast \\ X A_w^\ast \end{bmatrix} A_w^{*2} \\
 \qquad \vdots \\
  \begin{bmatrix} B_w^\ast \\ X A_w^{\ast}  \end{bmatrix} A_w^{*k-1}\\
\end{array}\right].
\end{equation}
\end{lemma}

\bpr
In the sequel it is convenient to denote
\[
 \Omega_w(k-1) =
\begin{bmatrix} B_w & A_w B_w & \cdots & A_w^{k-1} B_w\end{bmatrix}
 :\bigoplus_{k=0}^{k-1}\mathcal{E} \rightarrow \mathcal{X}_w.
\]
Recall that $B_w B_w^* = I - A_w A_w^*$. Using this
we have
\begin{align*}
\Omega_w(k-1)\Omega_w(k-1)^*
&= \sum_{j=0}^{k-1} A_w^j B_w B_w^* A_w^{*j}
=  \sum_{j=0}^{k-1} A_w^j (I-A_w A_w^*) A_w^{*j}\\
&= \sum_{j=0}^{k-1} A_w^j  A_w^{*j} - \sum_{j=1}^{k} A_w^j  A_w^{*j}
= I-A_w^k A_w^{*k}.
\end{align*}
In other words,
\[
\Omega_w(k-1) \Omega_w(k-1)^* = I-A_w^k A_w^{*k}.
\]
This, together with $ Q = I - X^* X$, yields
\begin{align*}
I - A_w^k Q A_w^{*k} & = I - A_w^k (I - X^* X) A_w^{*k} \\
&=( I - A_w^k A_w^{*k} ) + A_w^k X^* X A_w^{*k}\\
&= \Omega_w(k-1)\Omega_w(k-1)^*  +A_w^k X^* X A_w^{*k}.
\end{align*}
This readily implies that there exists a unitary operator
$\Xi$ from the range of $I - A_w^k Q A_w^{*k}$ onto the range of
$\begin{bmatrix}
   \Omega_w(k-1) & A_w^k X^* \\
 \end{bmatrix}^*$ such that
 \begin{equation}\label{rootk}
 \Xi \big(I - A_w^k Q A_w^{*k}\big)^{\frac{1}{2}} =
 \begin{bmatrix}
   \Omega_w(k-1)^* \\
   X A_w^{*k} \\
 \end{bmatrix}.
 \end{equation}
 In particular,
 \begin{equation}\label{kerkjj}
 \kr\big(I - A_w^k Q A_w^{*k}\big)  =
 \kr \begin{bmatrix}
   \Omega_w(k-1)^* \\
   X A_w^{*k} \\
 \end{bmatrix}.
 \end{equation}
 By exploiting this fact, we will derive \eqref{kerk}.
In fact, to complete the proof, it remains to show that
if $\Omega_w(k-1)^* x= 0$ and $X A_w^{*k} x=0$, then
$X A_w^{*j} x=0$ for $j <k$.

Recall that $X = \Upsilon_v \Upsilon_w^*$.
Let $S$ be the unilateral shift on $\ell_+^2(\mathcal{E})$.
Notice that $S^{*} \Upsilon_w^* = \Upsilon_w^* A_w^*$.
For any positive integer $j>0$, we have
\begin{equation}\label{shiftc}
X A_w^{*j} = \Upsilon_v S^{*j} \Upsilon_w^* =
\Upsilon_v \begin{bmatrix}
             B_w^*A_w^{*j} \\
             B_w^*A_w^{*(j+1)} \\
             B_w^*A_w^{*(j+2)} \\
             \vdots \\
           \end{bmatrix}=
\begin{bmatrix}
  0 & \Upsilon_v \\
\end{bmatrix}\begin{bmatrix}
               \Omega_w(j-1)^* \\
               \Upsilon_w^*A_w^{*j} \\
             \end{bmatrix}.
\end{equation}
Assume that there exists a vector $x\in \mathcal{X}_w$ such that
$\Omega_w(k-1)^* x=0$ and $X A_w^{*k}x =0$.
In particular $ B_w^* A_w^{*\ell} x =0 $ for $ \ell =0, \ldots, k-1$.
For any $j<k$, we have
\begin{align*}
X A_{_w}^{*j}x &=
\begin{bmatrix}
  0 & \Upsilon_v \\
\end{bmatrix}\begin{bmatrix}
               \Omega_w(j-1)^* x\\
               \Upsilon_w^* A_w^{*j}x \\
             \end{bmatrix}\\
&=\begin{bmatrix}
  0 & \star & A_v^{k-j}\Upsilon_v \\
\end{bmatrix}\begin{bmatrix}
               \Omega_w(j-1)^*x \\
               \begin{bmatrix}
               B_w^* A_w^{*j}x \\
               \vdots\\
               B_w^* A_w^{*k-1}x\\
               \end{bmatrix}\\
               \Upsilon_w^*A_w^{*k}x \\
             \end{bmatrix}  \\
&= A_v^{k-j} \Upsilon_v \Upsilon_w^* A_w^{*k}x =
A_v^{k-j} X A_w^{*k}x =0.
\end{align*}
Here $\star$ represents an unspecified entry.
So, if $X A_w^{*k}x =0$ then $X A_w^*A_w^{*j}x =0$
for all $j=0,1,2,\cdots,k-1$.
This, together with \eqref{kerkjj} yields \eqref{kerk} and
completes the proof. \epr

\bigskip

We conclude that Theorem \ref{mainthm00} implies Theorem 3.3 in \cite{GKRb} and vice versa.


\setcounter{equation}{0}
\section{Appendix: A unitary lower triangular operator}
Let us begin with the following result, which can be viewed
as a special case of Lemma 2.1 Page 76 in \cite{FF}.

\begin{lemma}\label{lem-unitary} Let $L$ be an operator of the form:
\begin{equation}\label{Llower}
L = \begin{bmatrix}
      A & 0 \\
      B & C \\
    \end{bmatrix}: \begin{bmatrix}
      \mathcal{X} \\
      \mathcal{Y} \\
    \end{bmatrix}\rightarrow \begin{bmatrix}
      \mathcal{X} \\
      \mathcal{Y} \\
    \end{bmatrix}.
\end{equation}
Then $L$ is a unitary operator if and only if the following
three conditions hold:
\begin{enumerate}
  \item $A$ is a co-isometry on $\mathcal{X}$.
  \item $C$ is an isometry on $\mathcal{Y}$.
  \item  The operator $B$ admits a decomposition of the form
  $B = P_{\mathfrak{H}_2} V P_{\mathfrak{H}_1}$
where $V$ is a unitary operator mapping
$\mathfrak{H}_1 = \mathcal{X}\ominus \im( A^*)$
onto the subspace $\mathfrak{H}_2 = \mathcal{Y}\ominus \im(C)$.
\end{enumerate}
\end{lemma}

\noindent {\sc Proof.} For completeness a proof is given.
Assume that $L$ is a unitary operator. Then using $L^*L= I$ and $LL^* = I$,
we have
\begin{align*}
L^*L &= \begin{bmatrix}
      A^* & B^* \\
      0 & C^* \\
    \end{bmatrix}\begin{bmatrix}
      A & 0 \\
      B & C \\
    \end{bmatrix} = \begin{bmatrix}
      A^*A + B^*B  & B^*C \\
      C^*B & C^*C \\
    \end{bmatrix} = \begin{bmatrix}
      I & 0 \\
      0 & I\\
    \end{bmatrix}\\
LL^* &=
\begin{bmatrix}
      A & 0 \\
      B & C \\
    \end{bmatrix} \begin{bmatrix}
      A^* & B^* \\
      0 & C^* \\
    \end{bmatrix} = \begin{bmatrix}
      A A^*    & AB^* \\
      BA^* & BB^* + CC^*\\
    \end{bmatrix} = \begin{bmatrix}
      I & 0 \\
      0 & I\\
    \end{bmatrix}.
\end{align*}
Hence $C^*C = I$ and $AA^*=I$. In other words, $C$ is an isometry and
$A$ is a co-isometry. This readily implies that
$P_{\mathfrak{H}_1} = I- A^*A$ is the orthogonal projection onto the subspace
$\mathfrak{H}_1 = \mathcal{X}\ominus \im(A^*)=\kr(A)$, and
$P_{\mathfrak{H}_2} = I- CC^*$ is the orthogonal projection onto the subspace
$\mathfrak{H}_2 = \mathcal{Y}\ominus \im( C) = \kr(C^*)$.

Since $B^*B = I- A^*A = P_{\mathfrak{H}_1}$, it follows that
$\im(B^*) =  \mathfrak{H}_1$. Moreover, there exists
a unitary operator $V$ mapping $\mathfrak{H}_1$ onto $\im(B)$ such that
$B = V P_{\mathfrak{H}_1}$. Using
$B B^* = I- CC^* = P_{\mathfrak{H}_2}$, it follows that
$\im(B)  =  \mathfrak{H}_2$.
Hence $B = P_{\mathfrak{H}_2} B = P_{\mathfrak{H}_2}V P_{\mathfrak{H}_1}$
where $V$ is a unitary operator mapping  $\mathfrak{H}_1$ into $\mathfrak{H}_2$.
Therefore Parts 1 to 3 hold.

On the other hand, if Parts 1 to 3 hold, then a direct calculation
shows that $L^*L = I$ and $LL^* =I$, and thus, $L$ is a unitary operator.
This completes the proof. \epr

\medskip

Let  $\Theta(z)$
 be a bi-inner rational function in
$H^\infty(\mathcal{E},\mathcal{E})$ whose Taylor series expansion is
given by $\Theta(z) = \sum_0^\infty z^n \Theta_n$. Then the corresponding
Laurent operator $L_\Theta$ is unitary and admits a lower triangular
matrix representation of the form:
\begin{equation}\label{L-Theta}
L_\Theta = \begin{bmatrix}
              \begin{bmatrix}
                \vdots & \vdots & \vdots & \vdots \\
                \cdots  & \Theta_0 & 0 & 0 \\
                \cdots & \Theta_1 & \Theta_0 & 0 \\
                \cdots & \Theta_2 & \Theta_1 & \Theta_0 \\
              \end{bmatrix}
              & \mbox{{ \Huge    0}} \\[1cm]
              \begin{bmatrix}
                \cdots & \Theta_3 & \Theta_2 & \Theta_1 \\
                \cdots & \Theta_4 & \Theta_3 & \Theta_2 \\
                 \cdots & \Theta_5 & \Theta_4 & \Theta_3 \\
                \vdots & \cdots  & \vdots & \vdots \\
              \end{bmatrix}
              & \begin{bmatrix}
                   \Theta_0 & 0 & 0 & \cdots \\
                   \Theta_1 & \Theta_0 & 0 & \cdots \\
                   \Theta_2 & \Theta_1 & \Theta_0 & \cdots \\
                   \vdots &  \vdots & \vdots  & \vdots \\
                 \end{bmatrix}
               \\
           \end{bmatrix}:\begin{bmatrix}
                           \ell_-^2(\mathcal{E}) \\
                           \ell_+^2(\mathcal{E}) \\
                         \end{bmatrix}\rightarrow\begin{bmatrix}
                           \ell_-^2(\mathcal{E}) \\
                           \ell_+^2(\mathcal{E}) \\
                         \end{bmatrix}.
\end{equation}

Let $T_\Theta$ on $\ell_+^2(\mathcal{E})$  be the lower triangular
Toeplitz matrix determined by
$\Theta$. Let $\widetilde{\Theta}$ be the function in $H^\infty(\mathcal{E},\mathcal{E})$
determined by $\widetilde{\Theta}(z) = \Theta(\overline{z})^* = \sum_0^\infty z^n \Theta_n^*$
when $|z| \leq 1$.
Now let $J$ be the unitary operator from
$\ell_+^2(\mathcal{E})$ onto $\ell_-^2(\mathcal{E})$ defined by
\[
J \begin{bmatrix}
    f_0 & f_1 & f_3 & \cdots \\
  \end{bmatrix}^{tr} = \begin{bmatrix}
    \cdots & f_2 & f_1 & f_0 \\
  \end{bmatrix}^{tr}
\]
where $^{tr}$ denotes the transpose.
Then we have
\[
J^*  \begin{bmatrix}
                \vdots & \vdots & \vdots & \vdots \\
                \cdots  & \Theta_0 & 0 & 0 \\
                \cdots & \Theta_1 & \Theta_0 & 0 \\
                \cdots & \Theta_2 & \Theta_1 & \Theta_0 \\
              \end{bmatrix} J = T_{\widetilde{\Theta}}^*.
\]
Notice that $T_{\widetilde{\Theta}}^*$ is the upper triangular Toeplitz
matrix on $\ell_+^2(\mathcal{E})$ determined by $\{\Theta_n\}_0^\infty$.
Moreover, $H_\Theta$ is the Hankel  operator on $\ell_+^2(\mathcal{E})$ defined by 
\begin{equation}\label{defHankel}
H_\Theta =
\begin{bmatrix}
                \Theta_1 & \Theta_2 & \Theta_3 & \cdots \\
                \Theta_2 & \Theta_3 & \Theta_4 & \cdots \\
                 \Theta_3 & \Theta_4 & \Theta_5 & \cdots \\
                \vdots & \vdots  & \vdots & \vdots \\
              \end{bmatrix} =
\begin{bmatrix}
                \cdots & \Theta_3 & \Theta_2 & \Theta_1 \\
                \cdots & \Theta_4 & \Theta_3 & \Theta_2 \\
                 \cdots & \Theta_5 & \Theta_4 & \Theta_3 \\
                \vdots & \vdots  & \vdots & \vdots \\
              \end{bmatrix} J.
\end{equation}
Using this in \eqref{L-Theta}, we see that
\begin{equation}\label{L-theta2}
 \begin{bmatrix}
              J^* & 0 \\
              0 & I \\
            \end{bmatrix} L_\Theta
            \begin{bmatrix}
              J & 0 \\
              0 & I \\
            \end{bmatrix} = \begin{bmatrix}
                              T_{\widetilde{\Theta}}^* & 0 \\
                              H_\Theta & T_\Theta \\
                            \end{bmatrix}:\begin{bmatrix}
                              \ell_+^2(\mathcal{E}) \\
                              \ell_+^2(\mathcal{E}) \\
                            \end{bmatrix} \rightarrow \begin{bmatrix}
                              \ell_+^2(\mathcal{E}) \\
                              \ell_+^2(\mathcal{E}) \\
                            \end{bmatrix}
\end{equation}
is a lower triangular unitary operator.
Observe  that both $T_\Theta$ and $T_{\widetilde{\Theta}}$ are isometries.
Now let $\mathfrak{H}(\Theta)$ and $\mathfrak{H}(\widetilde{\Theta})$
be the subspaces of $\ell_+^2(\mathcal{E})$ defined by
\begin{equation}\label{defHHH}
\mathfrak{H}(\Theta) =
\ell_+^2(\mathcal{E}) \ominus \im\left( T_\Theta\right)  \quad \mbox{and}\quad
 \mathfrak{H}(\widetilde{\Theta}) =
 \ell_+^2(\mathcal{E}) \ominus \im\left(T_{\widetilde{\Theta}}\right).
 \end{equation}
 By applying Lemma \ref{lem-unitary} to the lower triangular unitary
 matrix in \eqref{L-theta2}, we see that the Hankel
 operator $H_\Theta$ admits a decomposition of the form
 \begin{equation}\label{Hankeldec}
 H_\Theta = P_{\mathfrak{H}(\Theta)}VP_{\mathfrak{H}(\widetilde{\Theta})}
 = \big(I- T_\Theta T_\Theta^*\big)V
 \big(I- T_{\widetilde{\Theta}} T_{\widetilde{\Theta}}^*\big)
 \end{equation}
 where $V$ is a unitary operator from $\mathfrak{H}(\widetilde{\Theta})$
 onto $\mathfrak{H}(\Theta)$.

 Because the $2\times 2$ block  matrix in \eqref{L-theta2}
 is unitary, $\begin{bmatrix}
                H_\Theta & T_\Theta \\
              \end{bmatrix}$ is a co-isometry.
 Therefore
 \begin{equation}\label{coisoHT}
  T_\Theta T_\Theta^* + H_\Theta H_\Theta^* = I.
 \end{equation}
 In particular, $H_\Theta H_\Theta^*$ is the orthogonal
 projection onto $\mathfrak{H}(\Theta)$, that is,
 \begin{equation}\label{projHH}
 H_\Theta H_\Theta^* = I -  T_\Theta T_\Theta^* = P_{\mathfrak{H}(\Theta)}.
 \end{equation}

 Recall  that $H_\Theta^* = H_{\widetilde{\Theta}}$.
 Therefore the Hankel
 operator $H_{\widetilde{\Theta}}$ admits a decomposition of the form
 \begin{equation}  
 H_{\widetilde{\Theta}} = H_\Theta^*
  = P_{\mathfrak{H}(\widetilde{\Theta})}V^*P_{\mathfrak{H}(\Theta)}
 =  \big(I- T_{\widetilde{\Theta}} T_{\widetilde{\Theta}}^*\big)V^*
\big(I- T_\Theta T_\Theta^*\big)
 \end{equation}
 where $V$ is a unitary operator from $\mathfrak{H}(\widetilde{\Theta})$
 onto $\mathfrak{H}(\Theta)$. Finally,
 $H_{\widetilde{\Theta}} H_{\widetilde{\Theta}}^*$ is the orthogonal
 projection onto $\mathfrak{H}(\widetilde{\Theta})$, that is,
 \begin{equation}\label{projHH2}
 H_{\widetilde{\Theta}} H_{\widetilde{\Theta}}^* = I -
 T_{\widetilde{\Theta}} T_{\widetilde{\Theta}}^* = P_{\mathfrak{H}(\widetilde{\Theta})}.
 \end{equation}

\subsection*{Acknowledgements}

This work is based on research supported in part by the National Research Foundation of South Africa (NRF), (Grant Number 145688).

\end{document}